\documentstyle[11pt,amssymb,twoside]{article}
\def\TPP{{\rm T}\PP^2}

\def\ext{\mathop{\rm Ext}\nolimits}

\def\pic{\mathop{\rm Pic}\nolimits}
\def\co{{\cal O}}
\def\ch{\mathop{\rm ch}\nolimits}
\def\lll#1,#2{{\cal L}_{#1\div #2}}
\def\LL#1,#2{{\cal L}_{#1,#2}}
\newsavebox{\hdts}
\newsavebox{\ldts}
\newsavebox{\rdts}
\newsavebox{\pnt}
\sbox{\pnt}{\hbox to 0.2em{.\hss}}
\def\pline#1,#2{\raisebox{#1}[0em][0em]{\hbox to #2{\leaders\copy\pnt\hfil}}}
\newlength{\plw}

\def\upl#1{\raisebox{-0.1em}{\vbox{\hbox to \plw{\hfil$#1$\hfil}%
\hbox{\pline {-%
0.35%
em},{\plw}}}}}
\def\shdown#1{\raisebox{-0.2em}[0.8em][0em]{$#1$}}
\def\shup#1{\raisebox{0.2em}[0.4em][0.4em]{$#1$}}
\def\No{No.\ }
\def\x#1{{\bf(#1)}}
\def\xx#1{_{\mbox{\scriptsize\bf(#1)}}}

\def\le{\leqslant}
\def\ge{\geqslant}
\def\Ext{\mbox{\rm Ext}}
\def\Hom{\mbox{\rm Hom}}
\def\PP{{\Bbb P}}
\def\ZZ{{\Bbb Z}}

\def\EE{{\cal E}}
\def\FF{{\cal F}}
\def\GG{{\cal G}}
\def\OO{{\cal O}}
\def\ll{\langle}
\def\rr{\rangle}
\def\bigdot{{\mbox{\boldmath$\cdot$}}}
\def\chiminus{\chi\raisebox{-.25ex}{$_-$}}
\def\longto{\longrightarrow}
\def\ds{\displaystyle}
\newlength{\myskip}
\myskip=\abovedisplayskip
\author{B.~V.~Karpov$^1$ \and D.~Yu.~Nogin$^2$}
\title{Three-block exceptional collections over Del Pezzo surfaces}
\date{}

\begin{document}

\maketitle
\thispagestyle{empty}
\vskip -1.\baselineskip
{\renewcommand\thefootnote{}
\footnotetext{$^1$Partly supported by the International Science Foundation
(ISF grant No.~MKU\,300) and the Russian Fundamental Research Foundation
(project \No 95--01--00840).\\
\hspace*{18pt}$^2$Partly supported by the International Science Foundation
(ISF grant No.~MPN\,000) and the Russian Fundamental Research Foundation
(project \No 93--012--458).}}

\begin{abstract}
We study complete exceptional collections of coherent
sheaves over Del Pezzo surfaces, which consist of three
blocks such that inside each block all Ext groups between
the sheaves are zero. We show that the ranks of
all sheaves in such a block are the same and the three ranks corresponding to
a complete 3-block exceptional collection satisfy a Markov-type Diophantine
equation that is quadratic in each variable. For each Del Pezzo surface,
there is a finite number of these equations; the complete list is given. The
3-string braid group acts by mutations on the set of complete 3-block
exceptional collections. We describe this action. In particular, any orbit
contains a 3-block collection with the sum of ranks that is minimal for the
solutions of the corresponding Markov-type equation, and the orbits can be
obtained from each other via tensoring by an invertible sheaf and with the
action of the Weyl group. This allows us to compute the number of orbits up
to twisting.
\end{abstract}

\bigskip
\begin{center}{\large\bf Introduction}\end{center}

\medskip
Recall that a sheaf $E$ is called {\it exceptional\/} if
$\Hom(E,E)\cong{\Bbb C}$ and $\Ext^i(E,E)=0$ for $i>0$. An ordered collection
of sheaves $(E_1,\dots,E_\alpha)$ is called {\it exceptional\/} if all
$E_j$'s are exceptional and $\Ext^i(E_k,E_j)=0$ $\forall i$, $\forall k>j$.

The theory of exceptional sheaves is being developed for ten years. For the
first time, exceptional vector bundles over $\PP^2$ appeared in the paper
by J.-M. Drezet and J. Le Potier in 1985 [6]. In this paper, with the help of
discrete parameters $(r,c_1,c_2)$ of exceptional vector bundles, the boundary
of the set of these parameters for semistable sheaves over $\PP^2$ was
constructed (here, $r$ is the rank, $c_1,c_2$ are the Chern classes). In the
subsequent papers by Drezet [7--9], exceptional bundles were used in studying
the moduli space of semistable sheaves over $\PP^2$.

In another direction, the theory of exceptional sheaves was being developed
in Moscow by the participants of the seminar of Profs.\ A. N. Rudakov and A.
N. Tyurin. Originally, a general setting of the problem was to describe in a
reasonable way the set of exceptional sheaves over a given variety. Towards
this end, mutations of exceptional sheaves were used yielding new
exceptional sheaves. Also, the notion of a helix was introduced that ensured
the existence of mutations [3]. The first and most brilliant results were
obtained for $\PP^2$ [4]; namely, any exceptional sheaf is contained in some
helix and any helix can be obtained from a helix consisting of invertible
sheaves by a finite number of mutations. The possibility to obtain a given
exceptional collection from a certain canonical-form elementary one (e.g.,
for the case of $\PP^2$, from that consisting of invertible sheaves) by
mutations is called {\it constructivity}. The proof of the constructivity of
helices over $\PP^2$ is based on the following fact: The ranks $x,y,z$ of
three successive bundles in a helix obey the Markov equation
$$
x^2+y^2+z^2=3xyz;
\eqno(1)
$$
moreover, helix mutations exactly correspond to numerical mutations
of the equation solutions. As for the solutions in positive integers,
A.~A.~Markov [15] has shown that they can be reduced 
by mutations to the solution
$(1,1,1)$ with the minimum sum $x+y+z$.

This correspondence seems to be unique in some sense since for all known 
cases except $\PP^2$, there is no universal Diophantine equation on the ranks
of elements of a foundation of a helix. Therefore, to examine the
constructivity problem for other varieties, different methods were used.
Thus, the proof of the constructivity of helices and exceptional bundles over
$\PP^1\times\PP^1$ due to Rudakov [19] was based on geometric constructions
on the plane ${\rm Pic}\,(\PP^1\times\PP^1)\otimes{\Bbb Q}$, these
constructions being so subtle that they could hardly be employed for other
varieties. In the same paper, symmetric helices were studied, i.e., those
invariant under the involution which transposes the systems of generators
of $\PP^1\times\PP^1$. Of four bundles forming a foundation of a helix, two
are invariant under this involution and two others are mapped one onto
another. Let $x$ and $y$ be the ranks of the first two bundles and $z$ be the
rank of two others. Then, as is shown in [19], the triple $(x,y,z)$ is a
solution of the Diophantine equation
$$
x^2+y^2+2z^2=4xyz.
\eqno(2)
$$
All positive integer solutions of this equation can also be reduced by
mutations to $(1,1,1)$.

Using same methods, S. Yu.\ Zyuzina [10] proved the constructivity of
exceptional pairs over $\PP^1\times\PP^1$, i.e., the fact that any exceptional
pair is contained in some helix and, therefore, can be obtained by mutations
of invertible sheaves.

One of the authors devised the technique allowing one to prove
constructivity of not exceptional collections themselves, but their
images in the Grothendieck group K$_0$. This technique proved to be
especially fruitful for varieties with ${\rm rk\,K}_0=4$, namely, rational
ruled surfaces, $\PP^3$, 3-dimensional quadric, and Fano 3-folds
$V_5$ and $V_{22}$ (see [16, 17]).

Presently, the approach due to S. A. Kuleshov using rigid ($\Ext^1(F,F)=0$)
and superrigid ($\Ext^i(F,F)=0,\ i>0$) sheaves seems to be most effective.
Thus, over a Del Pezzo surface, any rigid sheaf is isomorphic to a direct
sum of exceptional sheaves. Basing on this, the constructivity of exceptional
collections of any length over a Del Pezzo surface is proved in [11]. The
extension of this result to the case of a surface with anticanonical class
free of basis components is presented in [12]. With the help of this
approach, Kuleshov obtained a new solution of the constructivity problem for
$\PP^2$ not employing the Markov equation (see [13]). This can open prospects
in investigating the case of $\PP^n$.

The aim of the present paper is to describe such classes of exceptional
collections over Del Pezzo surfaces, for which Markov-type equations exist
that play the same role as equations (1) and (2) play for $\PP^2$ and
$\PP^1\times\PP^1$ respectively. These are 3-block exceptional collections.

An exceptional collection $\EE=(E_1,\dots,E_\alpha)$ such that
$\Ext^i(E_j,E_k)=0$ $\forall i$, $\forall j\ne k$, is called a {\it block}.
In Sec.\ 1, we show that the ranks of all sheaves contained in a block are
equal; this number will be denoted by $r(\EE)$. A {\it 3-block collection\/}
is an exceptional collection
$$
(\EE,\FF,\GG)=(E_1,\dots,E_\alpha,F_1,\dots,F_\beta,G_1,\dots,G_\gamma)
$$
consisting of three blocks. The mutations of such collections preserving the
3-block structure are described in Sec.\ 2. These mutations define the action
of the 3-string braid group on the set of 3-block collections (see Sec.\
2.6).

The subject of our further study is not all 3-block collections but only
complete ones, i.e., those generating the derived category of coherent
sheaves over a surface under consideration. The main result of Sec.\ 3
is that the ranks $x=r(\EE)$, $y=r(\FF)$, $z=r(\GG)$ of a complete 3-block
collection $(\EE,\FF,\GG)$ obey the Markov-type equation
$$
\alpha x^2+\beta y^2+\gamma z^2=\sqrt{K^2\alpha\beta\gamma}\,xyz,
\eqno(3)
$$
where $K^2$ is the self-intersection index of the canonical class of a
surface. As is readily seen, the coefficients in the left-hand side are the
numbers of sheaves in the blocks $\EE$, $\FF$, $\GG$. The coefficient
$\sqrt{K^2\alpha\beta\gamma}$ in a right-hand side is an integer. This
imposes certain restrictions on $\alpha$, $\beta$, and
$\gamma$, thereby making the complete list of possible equations (3) for Del
Pezzo surfaces quite limited. The list is given in Sec.\ 3.5.
There, we also show that any positive integer solution of each equation
(3) can be obtained by mutations from the minimum solution, i.e., the
solution with the minimum sum $x+y+z$. This means that any 3-block collection
with given $\alpha$, $\beta$, and $\gamma$ can be obtained by mutations
preserving the 3-block structure from the collection with the same
$\alpha$, $\beta$, and $\gamma$ and minimum sum of ranks. Note that the
equations (1) and (2) are particular cases of (3).

In Sec.\ 4, we show that for each equation (3), a complete 3-block collection
with given $\alpha$, $\beta$, and $\gamma$ exists which corresponds to the
minimum solution of this equation. Moreover, any two such collections can be
obtained one from another by action of the Weyl group and tensoring by an
invertible sheaf (Sec.\ 5). This implies the following: Consider the orbits
of the action of the braid group on the set of complete 3-block collections.
If one joins into one class the orbits that differ by twisting only, the
number of classes corresponding to a fixed equation (3) will be finite. These
numbers are computed in Sec.\ 5.7.

\medskip
The authors would like to express their gratitude to A. L. Gorodentsev
and S. A. Kuleshov for the attention and interest to the subject.

\bigskip
\begin{center}{\large\bf 1. Preliminaries}\end{center}

\medskip
{\bf 1.1. Notations and agreements. }In what follows, $S$ is a Del Pezzo
surface over ${\Bbb C}$. By definition (see [14, Chapter~IV]), $S$ is a
birationally trivial smooth surface with ample anticanonical class
$(-K_S)$.

For a coherent sheaf $F$, we use the following discrete invariants:
the rank and Chern classes $r(F)$, $c_1(F)$, $c_2(F)$,
the slope $\mu(F)=\ds{c_1(F)\cdot(-K_S)\over r(F)}$,
$\nu(F)=\ds{c_1(F)\over r(F)}$ (a point in
${\rm Pic}\,(S)\otimes{\Bbb Q}$), the degree $d(F)=c_1(F)\cdot(-K_S)$.

For any coherent sheaves $E$ and $F$, denote by $\chi(E,F)$
the alternated sum
$$
\chi(E,F)=\sum_i(-1)^i\dim\Ext^i(E,F)
$$
which defines a bilinear form over K$_0(S)$. For torsion-free sheaves,
the relation
$$
\chi(E,F)=r(E)r(F)\left(\chi(\OO_S)+\frac{\mu(F)-\mu(E)}2+q(F)+q(E)-
{c_1(E)\cdot c_1(F)\over r(E)r(F)}\right)
\eqno(4)
$$
holds, which follows from the Riemann--Roch theorem.
Here, $q(F)=\frac{c_1^2(F)-2c_2(F)}{2r(F)}=\frac{\ch_2(F)}{r(F)}$.
Note that for Del Pezzo surfaces, $\chi({\cal O}_S)=1$.

The existence of a locally free resolvent for any coherent sheaf and
additivity of both sides of (4) imply the validity of this relation
for any two coherent sheaves. Removing the parentheses, we deduce
$$
\chi (E,F)=r(E)r(F)\chi({\cal O}_S)+\frac{r(E)d(F)}{2}-\frac{r(F)d(E)}{2}+
r(E)\ch_2(F)+r(F)\ch_2(E)-c_1(E)\cdot c_1(F).\quad
\eqno(5)
$$

The relation (4) immediately implies the expression for the skew-symmetric
part of the form $\chi$,
$$
\chiminus(E,F)\triangleq\chi(E,F)-\chi(F,E)=r(E)r(F)\big(\mu(F)-\mu(E)\big).
\eqno(6)
$$
If at least one of the ranks $r(E)$ and $r(F)$ is zero, one can use another
form of (6),
$$
\chiminus(E,F)=
\left|\begin{array}{cc} r(E) & r(F) \\ d(E) & d(F) \end{array}\right|.
\eqno(7)
$$

We denote the bounded derived category of coherent sheaves over $S$
by $D^b(S)$. Additive functions (in particular, $r$ and $d$) have natural
extensions over $D^b(S)$. Namely, if an object $A$ of the derived category
is represented by a complex $K^\bigdot$ with cohomology sheaves
$H^i(K^\bigdot)$ and $s$ is an additive function, then
$$
s(A)=\sum_i(-1)^is(K^i)=\sum_i(-1)^is(H^i(K^\bigdot)).
$$

\medskip
{\bf 1.2. Lemma. } {\it Let\/ $0\to E\to F\to G\to0$ be an exact sequence of
sheaves over\/ $S$. Then\/ $\chiminus(E,F)=\chiminus(F,G)=\chiminus(E,G)$.}

\smallskip
The statement of the lemma immediately follows from $(7)$.

\medskip
{\bf 1.3. Exceptional sheaves. }Recall that a sheaf $E$ is called {\it
exceptional\/} if $\Hom(E,E)\cong {\Bbb C}$, $\ext^i(E,E)=0$, $i>0$.

In [11], it is proved that an exceptional sheaf over a Del Pezzo surface
either is locally free, or is a torsion sheaf of the form $\OO_\ell(m)$, where
$m\in{\Bbb Z}$ and $\ell$ is an exceptional curve (or $(-1)$-curve),
i.e., an irreducible rational curve with $\ell^2=\ell\cdot K_S=-1$.

An exceptional vector bundle is Mumford--Takemoto stable with respect to
$(-K_S)$, and therefore, is uniquely determined up to an isomorphism by its
point $\nu$ (see~[5]). On the other hand, computation shows that
$\ch_2(\OO_\ell(m))=\frac{1}{2}-m$. Thus, the statement below is valid.

\smallskip
{\bf Proposition. }{\it An exceptional sheaf over a Del Pezzo surface\/ $S$
is uniquely determined up to an isomorphism by its image in\/ K$_0(S)$}.

\smallskip
The notion of exceptionality is naturally extended to the derived category.
An object $A\in D^b(S)$ is called {\it exceptional\/} if
$\Hom^0_{D^b(S)}(A,A)\cong {\Bbb C}$, $\Hom^i_{D^b(S)}(A,A)=0$, $i\ne 0$.
It is known (see [11]) that an object of a derived category over a Del Pezzo
surface is exceptional if and only if it is isomorphic to
$\delta E[i]$, where $E$ is an exceptional sheaf, $\delta$ denotes the
canonical embedding of the ground category in the derived one, $[i]$ denotes
the translation in $D^b(S)$.

\smallskip
Define the slope of an exceptional torsion sheaf as
$$
\mu(\OO_\ell(m))=+\infty.
$$

\medskip
{\bf 1.4. Lemma. }{\it Let\/ $E$ and\/ $F$ be exceptional sheaves over\/ $S$.
Then}
$$
{\rm sgn}\,\chiminus(E,F)={\rm sgn}\,(\mu(F)-\mu(E)).
$$

\smallskip
{\sc Proof}. For $r(E)r(F)\ne0$, the statement immediately follows from
$(6)$. If $r(E)=0$, i.e., $E=\OO_\ell(m)$, we have $d(E)=-K\cdot\ell=1$,
and $(7)$ implies $\chiminus(E,F)=-r(F)$. Similarly, for $r(F)=0$, we have
$\chiminus(E,F)=r(E)$. Clearly, this implies the desired equality.

\medskip
{\bf 1.5. Exceptional collections. }Recall that an ordered collection of
sheaves $(E_1,\dots,E_n)$ is called {\it exceptional\/} if all $E_j$'s are
exceptional sheaves, and for $1\le j<k\le n$, $\Ext^i(E_k,E_j)=0$ for any
$i$. Similarly, an exceptional collection of objects of $D^b(S)$ is defined.

An exceptional collection of sheaves such that for $j\ne k$,
$\Ext^i(E_j,E_k)=0$ for any $i$ is called a {\it block}. Evidently, a
collection obtained from a block by a permutation of its elements is also a
block.

\smallskip
{\bf Definition. }An {\it $m$-block collection\/} is an exceptional
collection
$$
(\EE_1,\EE_2,\dots,\EE_m)=
(E_{11},\dots,E_{1\alpha_1},E_{21},\dots,E_{2\alpha_2},\dots,
E_{m1},\dots,E_{m\alpha_m})
$$
such that all its subcollections $\EE_i=(E_{i1},\dots,E_{i\alpha_i})$ are
blocks. The {\it type\/} and {\it structure\/} of such a collection are,
respectively, the ordered and unordered collections of numbers
$(\alpha_1,\alpha_2,\dots ,\alpha_m)$. We will sometimes call $\alpha_i$ the
{\it length\/} of $\EE_i$.

\medskip
{\bf 1.6. Proposition. }{\it For a block\/ $(E_1,\dots,E_\alpha)$},

(a) $r(E_1)=\dots=r(E_\alpha);$

(b) $d(E_1)=\dots=d(E_\alpha);$

\noindent ({\it these invariants are denoted by\/ $r(\EE)$ $($the\/
{\rm rank} of the block$)$ and\/ $d(\EE)$ respectively.})

(c) {\it if\/ $r(\EE)=0$, then\/
$\EE=\left({\cal O}_{\ell_1}(D),\dots ,{\cal O}_{\ell_\alpha}(D)\right)$,
where\/ $\ell_i$'s are pairwise nonintersecting exceptional curves;}

(d) {\it the divisors $c_{ij}=c_1(E_i)-c_1(E_j)$, $i\ne j$, satisfy the
relations\/ $c_{ij}^2=-2$ and\/ $c_{ij}\cdot K_S=0$}.

\smallskip
{\sc Proof}. The definition immediately implies that
$\chiminus(E_j,E_k)=0$ $\forall j\ne k$.

If $r(E_j)=0$, then $E_j\cong\OO_\ell(m)$, and $c_1(E_j)=\ell$ is a
$(-1)$-curve, whence $d(E_j)=1$. Therefore, by (7), $E_k$ cannot be locally
free and thus has the same invariants $r=0$ and $d=1$.

Next, assume that $E_j$ is locally free. Then $E_k$ also is, and by (6),
$$
{d(E_j)\over r(E_j)}={d(E_k)\over r(E_k)}.
$$
The restriction of the exceptional bundle $E$ to a smooth elliptic curve
from the linear system $|{-K_S}|$ is a simple bundle of degree
$d(E)$ (see, e.g., [11]), and the rank and degree of a simple bundle over an
elliptic curve are coprime [1]. Hence, the validity of (a) and (b) follows.

Let us prove (c). We have
$\EE=({\cal O}_{\ell_1}(m_1),\dots ,{\cal O}_{\ell_\alpha}(m_\alpha))$.
The formula (5) under $r(E)=r(F)=0$ yields $\chi (E,F)=-c_1(E)\cdot c_1(F)$.
Hence, $\ell_i\cdot\ell_j=0$, $i\ne j$. Finally,
${\cal O}_{\ell_i}(m_i)\cong{\cal O}_{\ell_i}(D)$, where
$D=-\sum\limits_{i=1}^\alpha m_i\ell_i$.

To prove (d), consider two cases.

(i) $r(\EE)=0$. By (c), we have $c_{ij}=\ell_i-\ell_j$, where
$\ell_i^2=\ell_j^2=\ell_i\cdot K_S=\ell_j\cdot K_S=-1$ and
$\ell_i\cdot\ell_j=0$. Now the desired statement is verified by direct
computations.

(ii) $r(\EE)=r\ne 0$. Then the equality $\chi (E_i,E_i)=1$ together
with (4) implies
$$
q(E_i)=\frac{1}{2}\left(\frac{1}{r^2}+\frac{c_1(E_i)^2}{r^2}-1\right)\ ,
\quad i=1,\dots ,\alpha.
$$
It follows from (a) and (b) that $\mu(E_i)=\mu(E_j)$. Then
$$
0=\chi (E_i,E_j)=r^2\left(1+q(E_i)+q(E_j)-\frac{c_1(E_i)\cdot
c_1(E_j)}{r^2}\right)=1+\frac{1}{2}\left(c_1(E_i)-c_1(E_j)\right)^2.
$$
Therefore, $c_{ij}^2=-2$. The equality $c_{ij}\cdot K_S=0$ is a direct
consequence of (b). The proposition is proved.

\smallskip
{\bf Corollary. }{\it Let\/
$(\EE,\FF)=(E_1,\dots,E_\alpha,F_1,\dots,F_\beta)$ be a two-block exceptional
collection. Then\/ $\chi(E_j,F_k)=\chiminus(E_j,F_k)$ does not depend on\/
$j\in\{1,\dots,\alpha\}$ and\/ $k\in\{1,\dots,\beta\}$.}

\smallskip
We denote this quantity by $\chi(\EE,\FF)$.

\medskip
{\bf 1.7. Mutations} (see [5]). Let $(A,B)$ be an exceptional pair of objects
of $D^b(S)$. Consider an object $L_A^DB$ which completes the canonical
morphism ${\rm RHom}\,(A,B)\otimes A\stackrel{\rm can}\longto B$ to a
distinguished triangle
$$
L_A^DB[-1]\longto {\rm RHom}\,(A,B)\otimes A\stackrel{\rm can}\longto
B\longto L_A^DB.
\eqno(8)
$$
It is known [5] that $(L_A^DB,A)$ is an exceptional pair in $D^b(S)$.

\smallskip
{\bf Definition. }The {\it left mutation\/} is the mapping
$(A,B)\mapsto(L_A^DB,A)$ of the set of exceptional pairs of objects of
$D^b(S)$ onto itself.

\smallskip
The left mutation of a pair $(A,B)$ is the pair $(L_A^DB,A)$. The
object $L_A^DB$ is referred to as the result of a mutation, result of a
shift, or just a (left) shift of $B$ over $A$.

Next, let $(E,F)$ be an exceptional pair of sheaves over $S$. It is known [11,
2.11] that among the spaces $\Ext^i(E,F)$, either only
\Hom, or only $\Ext^1$ can be nonzero. By the definition of an exceptional
pair, $\Ext^i(F,E)=0\ \forall i$, whence
$$
\chi(E,F)=\chiminus(E,F).
$$
Using Lemma 1.4, we obtain the known [5] classification of exceptional pairs
in terms of the slopes. Namely, any pair $(E,F)$ has one of the
following types:
\begin{eqnarray*} \mbox{\rm hom-pair:}\quad & \Hom\,(E,F)\ne0,\
\Ext^i(E,F)=0,\ i=1,2, &
\Longleftrightarrow \; \mu(E)<\mu(F); \\ \mbox{\rm ext-pair:}\quad &
\Ext^1(E,F)\ne0,\ \Ext^i(E,F)=0,\ i=0,2, & \Longleftrightarrow \;
\mu(E)>\mu(F); \\ \mbox{\rm zero-pair:}\quad & \Ext^i(E,F)=0\ \forall i &
\Longleftrightarrow \; \mu(E)=\mu(F).
\end{eqnarray*}

Thus, a distinguished triangle (8) for a pair $(A,B)=(\delta E,\delta F)$
is always reduced to one of the following exact triples:
\begin{eqnarray*}
0\longto L_EF\longto\Hom\,(E,F)\otimes E\longto F\longto0 & \quad &
\hbox{(division);} \\
0\longto\Hom\,(E,F)\otimes E\longto F\longto L_EF\longto0 & \quad &
\hbox{(recoil);} \\
0\longto F\longto L_EF\longto\Ext^1(E,F)\otimes E\longto0 & \quad &
\hbox{(extension).}
\end{eqnarray*}

Since the pair $(L_{\delta E}^D\delta F,\delta E)$ of objects of $D^b(S)$
is exceptional, the same is true for the pair of sheaves $(L_EF,E)$.

\smallskip
{\bf Definition. }The {\it left sheaf mutation\/} is the mapping
$(E,F)\mapsto(L_EF,E)$ of the set of exceptional pairs of sheaves onto itself.
Three types of sheaf mutations are distinguished, namely, {\it division},
{\it recoil}, and {\it extension}, depending on which of the exact triples
given above takes place.

\smallskip
Dually, right mutations are defined; if $(A,B)$ is an exceptional pair in
$D^b(S)$, then the object $R_B^DA$ completes the canonical morphism
$A\to{\rm RHom}^*(A,B)\otimes B$ to the distinguished triangle
$$
R_B^DA\longto A\longto{\rm RHom}^*(A,B)\otimes B\longto R_B^DA[1].
$$
The mapping $(A,B)\mapsto(B,R_B^DA)$ is called the right mutation. The pair
$(B,R_B^DA)$ is also exceptional.

For an exceptional pair of sheaves $(E,F)$, the distinguished triangle above
is reduced to one of the following triples:
\begin{eqnarray*}
0\longto E\longto\Hom\,(E,F)^*\otimes F\longto R_FE\longto0 & \quad &
\hbox{(division);} \\
0\longto R_FE\longto E\longto\Hom\,(E,F)^*\otimes F\longto0 & \quad &
\hbox{(recoil);} \\
0\longto\Ext^1(E,F)^*\otimes E\longto R_FE\longto E\longto0 & \quad &
\hbox{(extension).}
\end{eqnarray*}

The right sheaf mutation is the mapping $(E,F)\mapsto(F,R_FE)$ of the set of
exceptional pairs onto itself.

The following obvious statement establishes the relation between mutations in
the derived category and sheaf mutations.

\medskip
{\bf 1.8. Proposition. }{\it If a sheaf mutation of a pair\/ $(E,F)$ is
either recoil or extension, then\/ $L_{\delta E}^D\delta F=\delta L_EF$ and\/
$R_{\delta F}^D\delta E=\delta R_FE$. For the case of division, $L_{\delta
E}^D\delta F=\delta L_EF[1]$ and\/ $R_{\delta F}^D\delta E=\delta R_FE[-1]$.}

\smallskip
Further on, we identify the category of coherent sheaves over
$S$ with its image in $D^b(S)$ under the canonical embedding. The symbol
$\delta$ is usually omitted. Under this agreement, the latter proposition
reads that a sheaf mutation other than a division coincides with the mutation
in $D^b(S)$, and for a division-type mutation, the shift by $\pm1$ in
$D^b(S)$ occurs.

\medskip
{\bf 1.9. Proposition. }{\it The type of a left sheaf mutation is
described in terms of slopes as follows:}
\begin{eqnarray*}
\hbox{(division)} & \Longleftrightarrow & \mu(L_EF)<\mu(E)<\mu(F); \\
\hbox{(recoil)} & \Longleftrightarrow & \mu(E)\le\mu(F)\le\mu(L_EF); \\
\hbox{(extension)} & \Longleftrightarrow & \mu(F)\le\mu(L_EF)\le\mu(E).
\end{eqnarray*}

\smallskip
{\sc Proof}. For $\mu(E)=\mu(F)$, i.e., $\Ext^i(E,F)=0\ \forall i$,
the left mutation of the pair $(E,F)$ can be regarded as either trivial
recoil or trivial extension and is just the transposition of $E$ and $F$.
Lemmas 1.2 and 1.4 yield $\mu(L_EF)=\mu(E)=\mu(F)$.

Now, let $\mu(E)\ne\mu(F)$. If $\mu(E)<\mu(F)$, i.e., $(E,F)$ is a
hom-pair, then the left mutation of this pair is either division or recoil.
For the case of division, the correspondent exact triple shows that the pair
$(L_EF,E)$ is also of hom type, whence the desired inequality follows. For a
recoil,
$$
\chiminus(F,L_EF)\stackrel{\rm 1.2}=\chiminus(E\otimes\Hom\,(E,F),F)=
        \dim\Hom(E,F)\cdot\chiminus(E,F)>0,
$$
and it remains to apply Lemma 1.4.

Finally, if $\mu(E)>\mu(F)$, i.e., $(E,F)$ is an ext-pair, then the left
mutation of this pair is an extension. The correspondent exact triple
shows that $(L_EF,E)$ is a hom-pair, and
$$
\chiminus(F,L_EF)=\chiminus(F,\Ext^1(E,F)\otimes E)=
        \dim\Ext^1(E,F)\cdot\chiminus(F,E)>0,
$$
which completes the proof.

\medskip
{\bf 1.10. }Let $\tau=(E_1,\dots,E_n)$ be an exceptional collection of sheaves
(or objects of $D^b(S)$). It is known [2, 5], that the mappings
\begin{eqnarray*}
\tau&\longmapsto&(E_1,\dots,E_{i-1},L_{E_i}^{(D)}E_{i+1},E_i,E_{i+2},
\dots,E_n)\\[\abovedisplayskip]
\noalign{\noindent\mbox{and}}\\[\myskip]
\tau&\longmapsto&
(E_1,\dots,E_{i-1},E_{i+1},R_{E_{i+1}}^{(D)}E_i,E_{i+2},\dots,E_n),
\end{eqnarray*}
$i=1,\dots,n-1$, result in exceptional collections. These mappings are called
{\it mutations of exceptional collections.}

\medskip
{\bf 1.11. Helices. }Consider a bi-infinite extension of an exceptional
collection of sheaves $(E_1,\dots,E_n)$ defined recursively by
$$
E_{i+n}=R_{E_{i+n-1}}\dots
R_{E_{i+1}}E_i,\qquad E_{-i}=L_{E_{1-i}}\dots L_{E_{n-1-i}}E_{n-i},\qquad
i\ge1.
$$
The sequence $\{E_m\}_{m\in{\Bbb Z}}$ thus constructed is called a (sheaf)
{\it helix of period\/} $n$ if
$$
E_i=E_{n+i}\otimes K_S \quad \forall i\in{\Bbb Z}.
$$

A {\it foundation of a helix\/} is any its subcollection of the form
$(E_{m+1},\dots,E_{m+n})$, where $m\in{\Bbb Z}$.

\smallskip
{\bf Definition. }An exceptional collection $(E_1,\dots,E_n)$ is called
{\it complete\/} if it generates $D^b(S)$.

\smallskip
For collections over a variety with ample anticanonical class,
it is known [2] that an exceptional collection is complete if and only if it
is a foundation of a helix.

Below, we study complete 3-block exceptional collections over Del Pezzo
surfaces.

\bigskip
\begin{center}{\large\bf 2. Mutations of block collections}\end{center}

\medskip
{\bf 2.1. }Let $(\EE,\FF)=(E_1,\dots,E_\alpha,F_1,\dots,F_\beta)$ ---
be a two-block collection. Define the {\it shift of\/ $F_j$ over\/}
$\EE$ as
$$
L_\EE F_j=L_{E_1}\dots L_{E_\alpha}F_j,\quad j=1,\dots,\beta.
$$
According to [5, Sec.\ 4], $\Ext^i(L_\EE F_j,L_\EE F_k)=\Ext^i(F_j,F_k)=0\
\forall i,\ \forall j\ne k$. Hence, the collection of sheaves
$$
L_\EE\FF\triangleq(L_\EE F_1,\dots,L_\EE F_\beta)
$$
is a block; we will sometimes call it the {\it shift of\/ $\FF$ over\/
$\EE$}. Consider the sheaves 
$$
R_\FF E_i=R_{F_\beta}\dots R_{F_1}E_i
$$ 
and the
collection $R_\FF\EE=(R_\FF E_1,\dots,R_\FF E_\alpha)$ which is a block too.

One can easily see that the two-block collections $(L_\EE\FF,\EE)$ and
$(\FF,R_\FF\EE)$ are obtained from the initial collection $(\EE,\FF)$ by a
finite number of mutations in the sense of 1.10 and, therefore, are
exceptional.

We will also use the block $L^D_\EE\FF$ obtained as a result of
replacing sheaf mutations in the definition of $(L_\EE\FF,\EE)$ with mutations
in the derived category,
$$
L^D_\EE\FF\triangleq(L^D_\EE F_1,\dots,L^D_\EE F_\beta),
\qquad\mbox{where}\quad
L^D_\EE F_j=L^D_{E_1}\dots L^D_{E_\alpha}F_j.
$$
\smallskip
{\bf Definition. }The mappings $(\EE,\FF)\mapsto(L_\EE\FF,\EE)$ and
$(\EE,\FF)\mapsto(\FF,R_\FF\EE)$ of the set of two-block collections onto
itself are called, respectively, the {\it left\/} and {\it right mutations
of two-block collections}.

\medskip
{\bf 2.2. Proposition. }1. {\it The sheaf $L_\EE F_j$ is contained in one of
the following exact triples:}
\begin{eqnarray*}
0\longto L_\EE F_j\longto
\bigoplus_{i=1}^\alpha\Big(\Hom\,(E_i,F_j)\otimes E_i\Big)
\stackrel{\rm can}\longto F_j\longto0 & \quad & \hbox{(division),}\\
0\longto\bigoplus_{i=1}^\alpha\Big(\Hom\,(E_i,F_j)\otimes E_i\Big)
\stackrel{\rm can}\longto F_j\longto L_\EE F_j\longto0
& \quad & \hbox{(recoil),}\\
0\longto F_j\longto L_\EE F_j\longto
\bigoplus_{i=1}^\alpha\Big(\Ext^1(E_i,F_j)\otimes E_i\Big)
\longto0 & \quad & \hbox{(extension).}
\end{eqnarray*}

2. {\it These three cases are described in terms of the discrete invariants
as follows:}
\begin{eqnarray*}
\hbox{(division)} & \Longleftrightarrow & \alpha\chi(\EE,\FF)r(\EE)>r(\FF),\\
\hbox{(recoil)} & \Longleftrightarrow & \chi(\EE,\FF)\ge0\quad\hbox{\it and}
        \quad\alpha\chi(\EE,\FF)r(\EE)\le r(\FF),\\
\hbox{(extension)} & \Longleftrightarrow & \chi(\EE,\FF)\le0,
\end{eqnarray*}
{\it In particular, the type of the exact triple of item\/ $1$ does not
depend on\/} $j\in\{1,\dots,\beta\}$.

\smallskip
{\sc Proof}. Identify the category of coherent sheaves over $S$ with its
image under the canonical embedding in $D^b(S)$. Consider the direct sum of
canonical morphisms
$$
\bigoplus_{i=1}^\alpha\Big({\rm RHom}\,(E_i,F_j)\otimes E_i\Big)\longto F_j
$$
which can be completed to a distinguished triangle
$$
X[-1]\longto
\bigoplus_{i=1}^\alpha\Big({\rm RHom}\,(E_i,F_j)\otimes E_i\Big)\longto F_j
\longto X.
\eqno(9)
$$
Applying the functor $\mbox{RHom}\,(E_{i'},\bigdot)$ to it and taking into account
that $\mbox{RHom}(E_{i'},E_i)=0$ for $i\ne i'$, we obtain the long 
exact sequence
\begin{eqnarray*}
\dots\!\!&\!\!\longto\!\!&\!\!\Hom_{D^b(S)}^{k-1}(E_{i'},X)\longto
\Hom_{D^b(S)}^k(E_{i'},{\rm RHom}\,(E_{i'},F_j)\otimes E_{i'})
\stackrel{f_k}\longto\Hom_{D^b(S)}^k(E_{i'},F_j)\longto \\
\!\!&\!\!\longto\!\!&\!\!\Hom_{D^b(S)}^k(E_{i'},X)\longto\dots,
\end{eqnarray*}
where all $f_k$'s are canonical identities due to the exceptionality of
$E_{i'}$. Hence, all spaces\linebreak $\Hom_{D^b(S)}^k(E_{i'},X)$ are zero,
i.e.,
$$
{\rm RHom}\,(E_{i'},X)=0\quad\forall i'\in\{1,\dots,\alpha\}.
\eqno(10)
$$

Next, applying the functor $\mbox{RHom}(F_j,\bigdot)$ to (9) and taking into
account that $\mbox{RHom}(F_j,E_i)=0\ \forall i\in\{1,\dots,\alpha\}$, 
we obtain
the long exact sequence consisting of spaces $\Hom^k_{D^b(S)}(F_j,\bigdot)$,
whence we conclude that
$$
\Hom^0_{D^b(S)}(F_j,F_j)\cong\Hom^0(F_j,X)\cong{\Bbb C}\quad\mbox{and}\quad
\Hom^k_{D^b(S)}(F_j,X)=0,\ k\ne0.
\eqno(11)
$$

Finally, applying $\mbox{RHom}\,(\bigdot\,,X)$ to (9), we obtain
$\Hom^0_{D^b(S)}(X,X)\cong\Hom^0(F_j,X)\cong{\Bbb C}$ and
$\Hom^k_{D^b(S)}(X,X)=0,\ k\ne0$, i.e., $X$ is an exceptional object in
$D^b(S)$.

Let $T={\rm Tr}\,(E_1,\dots,E_\alpha,F_j)$ be a complete triangled
subcategory in $D^b(S)$ generated by the corresponding collection of objects;
it contains $T_\EE={\rm Tr}\,(E_1,\dots,E_\alpha)$. The distinguished
triangle (9) and the equality (10) mean that $X$ belongs to the intersection
of the right orthogonal subcategory $T_\EE^\perp$ (see [2, Sec.~3]) and $T$,
this intersection being generated by the exceptional object $L_\EE F_j$.
Since $X$ is exceptional, it is quasi-isomorphic to $L_\EE F_j[p]$, the
latter being a complex with $L_\EE F_j$ in the $p$th position and zeroes in
the others.

Therefore, any object in (9) has only one nonzero cohomology sheaf; hence,
this distinguished triangle is reduced to one of the exact triples from the
proposition statement. This proves item 1.

The assertions of item 2 are obvious due to Lemma 1.4, except only for the fact that
under $\alpha\chi(\EE,\FF)r(\EE)=r(\FF)$, the division cannot take place.
Let us show this. Indeed, if the contrary holds, $r(L_\EE \FF)=0$, and hence,
$r(\EE)=0$ since a torsion sheaf is never a subsheaf of a locally free one.
Hence, $r(\FF)=0$ too. Then, by (7), $\chi(\EE,\FF)=0$, i.e., the mutation
is a trivial recoil and $L_\EE F_j=F_j$. Note that in this case the
collection $(\EE,\FF)$ is actually a single block. The proposition is proved.

\smallskip
Dual reasoning easily proves the analogous facts for right mutations.

\medskip
{\bf 2.3. Proposition. }1. {\it The sheaf\/ $R_\FF E_i$ is contained in one
of the exact triples}
\begin{eqnarray*}
0\longto E_i\stackrel{\rm can}\longto
\bigoplus_{i=1}^\beta\Big(\Hom^*(E_i,F_j)\otimes F_j\Big)
\longto R_\FF E_i\longto0 & \quad & \hbox{(division),}\\
0\longto R_\FF E_i\longto E_i\stackrel{\rm can}\longto
\bigoplus_{i=1}^\beta\Big(\Hom^*(E_i,F_j)\otimes F_j\Big)\longto0
& \quad & \hbox{(recoil),}\\
0\longto\bigoplus_{i=1}^\beta\Big(\Ext^1(E_i,F_j)^*\otimes F_j\Big)
\longto R_\FF E_i\longto E_i\longto0 & \quad & \hbox{(extension).}
\end{eqnarray*}

2. {\it These three cases can be described in terms of the discrete
invariants as follows:}
\begin{eqnarray*}
\hbox{(division)} & \Longleftrightarrow & r(\EE)\le\beta\chi(\EE,\FF)r(\FF),\\
\hbox{(recoil)} & \Longleftrightarrow & r(\EE)>\beta\chi(\EE,\FF)r(\FF)\ge0,\\
\hbox{(extension)} & \Longleftrightarrow & \chi(\EE,\FF)\le0,
\end{eqnarray*}
{\it In particular, the type of the exact triple of item\/ $1$ does not
depend on\/} $i\in\{1,\dots,\alpha\}$.

\smallskip
Note that, for the simplest case where each of the blocks $\EE$ and
$\FF$ consists of one sheaf, exact triples of Propositions 2.2 and 2.3
coincide with those of Sec.\ 1.7.

\medskip
{\bf 2.4. }A mutation of a two-block collection $(\EE,\FF)$ is called a
{\it division\/}, {\it recoil\/}, or {\it extension\/} depending on 
the type of the corresponding
exact triples given in the latter two propositions.

Our next aim is to show that a non-division-type mutation of a two-block
collection coincides with the mutation in $D^b(S)$, and under a division,
a shift of grading by $\pm1$ occurs.

Below, precise statements for left mutations are presented; we leave the case
of right mutation to the reader.

For convenience, renumber a block $\EE$ in the reverse order,
$\EE=(E_\alpha,\dots,\EE_1)$; then $L_\EE F_j=L_{E_\alpha}\dots L_{E_1}F_j$.
Consider the sequence of sheaves
$$
L^0F_j=F_j,\qquad L^iF_j=L_{E_i}\dots L_{E_1}F_j,\quad i=1,\dots,\alpha.
$$

Note that $L^\alpha F_j=L_\EE F_j$ is the shift of $F_j$ over the block
$\EE$, and $L^iF_j$ is the shift of $F_j$ over $\EE_i=(E_i,\dots,E_1)$.
Denote by $L^i$ the left (sheaf) mutation of the pair $(E_i,L^{i-1}F_j)$.
Then the sheaf $L_\EE F_j$ is a result of applying the sequence of mutations
$L^1,\dots,L^\alpha$.

\medskip
{\bf 2.5. Proposition. }1. {\it If the left mutation of a two-block
collection\/ $(\EE,\FF)$ is not a division, then all mutations\/ $L^i,\
i=1,\dots,\alpha,$ are not divisions for any\/ $F_j\in\FF$.}

2. {\it If the left mutation of \/$(\EE,\FF)$ is a division, then exactly one
mutation in the sequence\/ $L^i$ is a division.}

\smallskip
{\sc Proof}. Let $\chi(\EE,\FF)>0$. Then, by 2.2, the left mutation of the
two-block collection $(\EE,\FF)$ is either a division or a recoil. For the
case of a division, we obtain by Lemma 1.2 that
$$
\begin{array}{ccc}
\chiminus
\Big(L_\EE F_j,\bigoplus\limits_{i=1}^\alpha\Hom(E_i,F_j)\otimes E_i\Big)\ &
= & \chiminus
\Big(\bigoplus\limits_{i=1}^\alpha\Hom(E_i,F_j)\otimes E_i,F_j\Big), \\
\mbox{\scriptsize $|\vspace{0.12ex}|$} &&
\mbox{\scriptsize $|\vspace{0.12ex}|$} \\[4pt]
\alpha\chi(\EE,\FF)\chi(L_\EE\FF,\EE) && \alpha\chi(\EE,\FF)^2
\end{array}
$$
whence $\chi(L_\EE\FF,\EE)>0$ and $\mu(L_\EE\FF)<\mu(\EE)$.

Consider the sequence $\mu_i=\mu(L^iF_j)$, where
$\mu_0=\mu(\FF)>\mu(\EE)$ and $\mu_\alpha=\mu(L_\EE\FF)<\mu(\EE)$. Let $p$ be
the least number where the change of the sign of $\mu_i-\mu(\EE)$ occurs,
i.e.,
$$
\mu_{p-1}-\mu(\EE)>0\qquad\hbox{and}\qquad\mu_p-\mu(\EE)<0.
$$
By Proposition 1.9, each mutation $L^i$, $1\le i\le p-1$, can be none other
than a recoil, and $L^p$ is a division. Then, again by Proposition 1.9,
all mutations $L^i$, $p+1\le i\le\alpha$, are extensions.

Similarly, for the case where the left block mutation of
$(\EE,\FF)$ is a recoil, we have $\mu(L_\EE\FF)>\mu(\FF)>\mu(\EE)$.
Proposition 1.9 implies that all $L^i$'s are recoils in this case.

Let now $\chi(\EE,\FF)<0$. Computation shows that
$\mu(\FF)<\mu(L_\EE\FF)<\mu(\EE)$, and by Proposition 1.9, all mutations
$L^i$ are extensions.

Finally, for $\chi(\EE,\FF)=0$, evidently, all $L^i$'s are trivial recoils
(or trivial extensions). This completes the proof of the proposition.

\medskip
{\bf Corollary. }{\it An object\/ $X$ in the distinguished triangle\/ $(9)$
coincides with\/ $L^D_\EE\FF_j$.}

\smallskip
{\sc Proof}. According to the proof of Proposition 2.2, $X$ is
quasi-isomorphic to $L^D_\EE\FF_j[p]$ for some $p$. Applying the preceding
proposition and 1.8, we obtain that $L^D_\EE\FF_j=L_\EE\FF_j$ if the left
mutation of the two-block collection $(\EE,\FF)$ is not a division,
and $L^D_\EE\FF_j=L_\EE\FF_j[1]$ otherwise. Examining the correspondence
between the distinguished triangle (9) and exact triples of Proposition 2.2,
we obtain the desired statement.

\medskip
{\bf 2.6. Action of the braid group. }Let $\tau=(\EE_1,\dots,\EE_m)$ be an
$m$\/-block collection. Define {\it left\/} and {\it right mutations of
$m$-block collections\/} as mappings of the set of $m$\/-block collections
onto itself,
\begin{eqnarray*}
L_i:\tau&\longmapsto &
(\EE_1,\dots,\EE_{i-1},L_{\EE_i}\EE_{i+1},\EE_i,\EE_{i+2},\dots,\EE_m)
\\[\abovedisplayskip]
\noalign{\noindent\mbox{and}}\\[\myskip]
R_i:\tau&\longmapsto&
(\EE_1,\dots,\EE_{i-1},\EE_{i+1},R_{\EE_{i+1}}\EE_i,\EE_{i+2},\dots,\EE_m),
\end{eqnarray*}
where $i=1,\dots,m-1$.
Note that the structure of an $m$\/-block collection defined in 1.5 is
preserved under mutations $L_i$ and $R_i$. Hence, any orbit under the action
of the braid group is contained in the set of $m$\/-block collections of a
correspondent structure. We say that the mutations of
$\tau$ in the sense of 1.10 do not preserve the 3-block structure.
The statement below is an analog of [2, 2.3] for block collections.

\smallskip
{\bf Proposition. }1. {\it Mutations\/ $R_i$ and\/ $L_i$ are inverse, i.e.,\/}
$R_i\circ L_i={\rm id}$.

2. {\it Right\/ $($and left\/$)$ mutations define the action of the
$m$-string braid group, i.e., the generating relations of the braid
group hold,}
$$
R_i \circ R_{i+1}\circ R_i=R_{i+1}\circ R_i\circ R_{i+1}, \qquad
L_i \circ L_{i+1}\circ L_i=L_{i+1}\circ L_i\circ L_{i+1}.
$$

\smallskip
{\sl Proof}. To prove item 1, it suffices to check that for a two-block
collection $(\EE,\FF)$, the equality $R_\EE(L_\EE\FF)=\FF$ holds. This easily
follows from the fact that left sheaf mutations are inverse to right ones.

To prove item 2, it suffices to check the following statement: If
$(\EE,\FF,\GG)$ is a 3-block collection, then
$$
R_\GG(R_\FF\EE)=R_{R_\GG\FF}R_\GG\EE.
$$
Indeed, according to [2], the block in the left-hand side of the latter
equality is the right shift of $\EE$ over the category ${\rm
Tr}\,(\FF,\GG)={\rm Tr}\,(\GG,R_\GG\FF)$ and does not depend on the choice of
a basis in it.

\medskip
{\bf 2.7. }In conclusion, note that one can use a ``matrix notation'' for
mutations of two-block collections. For example, a division-type left
mutation of $(\EE,\FF)$ corresponds to the sequence
$$
\left(\!\!
\begin{array}{c}0\\ \vdots\\ 0\end{array}
\!\!\right)
\!\!\to\!\!
\left(\!\!
\begin{array}{c}L_\EE F_1\\ \vdots\\ L_\EE F_\beta\end{array}
\!\!\right)
\!\!\to\!\!
\left[\!\left(\!\!
\begin{array}{ccc}
\Hom(E_1,F_1) & \dots & \Hom(E_\alpha,F_1) \\
\multicolumn{3}{c}{\dotfill} \\
\Hom(E_1,F_\beta) & \dots & \Hom(E_\alpha,F_\beta)
\end{array}
\!\!\right)\!\odot\!\left(\!\!
\begin{array}{c}E_1\\ \vdots\\ E_\alpha\end{array}
\!\!\right)\!\right]
\!\!\to\!\!
\left(\!\!
\begin{array}{c}F_1\\ \vdots\\ F_\beta\end{array}
\!\!\right)
\!\!\to\!\!
\left(\!\!
\begin{array}{c}0\\ \vdots\\ 0\end{array}
\!\!\right)
\eqno(12)
$$
where the result of the ``multiplication'' $\odot$ of the matrix
$\Big(\Hom(E_i,F_j)\Big)$ by the column of sheaves is the column of sheaves
with
$$
\Big(\!\!
\begin{array}{ccc}
\Hom(E_1,F_j) & \dots & \Hom(E_\alpha,F_j)
\end{array}\!\!
\Big)\!\odot\!\left(\!\!
\begin{array}{c}E_1\\ \vdots\\ E_\alpha\end{array}\!\!
\right)
=\bigoplus_{i=1}^\alpha
\Big(\Hom(E_i,F_j)\otimes E_i\Big)
$$
in the $j$th position.

It is quite natural to denote a matrix consisting of vector spaces
$\Hom(E_i,F_j)$ by $\Hom(\EE,\FF)$. Then the middle term in (12)
takes the form $\Hom(\EE,\FF)\odot\EE$, and exact sequences like (12)
that correspond to various types of mutations of $(\EE,\FF)$ can be obtained
from the sequences of Sec.~1.7 (which define mutations of an exceptional
pair of sheaves $(E,F)$) by replacing $E$, $F$, and the symbol of tensor
product with $\EE$, $\FF$, and $\odot$.

\bigskip
\begin{center}
{\large\bf 3. Markov-type equations for complete 3-block collections}
\end{center}

\medskip
In this section, we always assume
$$
(\EE,\FF,\GG)=(E_1,\dots,E_\alpha,F_1,\dots,F_\beta,G_1,\dots,G_\gamma)
$$
to be a complete 3-block collection of sheaves over a Del Pezzo surface
$S$ (in the sense of 1.5). The completeness is equivalent to the fact that
this ordered collection of sheaves is a foundation of a helix (see 1.11).

\medskip
{\bf 3.1. }Consider a $\Bbb Z$-module K$_0(S)$ with the bilinear form
$\ll x,y\rr=\chi(x,y)$. Let $\lambda:{\rm K}_0(S)\to{\rm K}_0(S)^*$ be
defined as $\lambda x=\ll\bigdot\,,x\rr$. For any additive functions
$s$ and $t$ on K$_0(S)$, define
$$
\ll s,t\rr\triangleq\ll\lambda^{-1}s,\lambda^{-1}t\rr.
\eqno(13)
$$

Consider the additive functions $r$ and $d$, where $r$ is the rank and
$d(U)=c_1(U)\cdot(-K_S)$. Then $\lambda^{-1}r=\OO_p$, the latter being a
structure sheaf of a point. Direct computations show that $\lambda^{-1}d$
lies in the linear span of $\OO_{-K_S}$ and $\OO_p$. Therefore,
\begin{eqnarray*}
\ll r,r\rr & = & \chi(\OO_p,\OO_p)\>=\>r(\OO_p)\>=\>0, \\
\ll r,d\rr & = & \ll\OO_p,\lambda^{-1}d\rr\>=\>d(\OO_p)\>=\>0, \\
\ll d,r\rr & = & r(\lambda^{-1}d)\>=\>0.
\end{eqnarray*}

Let $(e_1,\dots,e_n)$ be a semiorthogonal basis of K$_0(S)$ (i.e.,
$\ll e_i,e_i\rr=1$, $\ll e_j,e_i\rr=0$, $j>i$), $(e_n^\vee,\dots,e_1^\vee)$
be a {\it dual\/} semiorthogonal basis, i.e., such that
$\ll e_i,e_j^\vee\rr=\delta_{ij}$. Then, for additive functions
$s$ and $t$, one has
$$
\ll s,t\rr=\sum_{i=1}^ns(e_i)t(e_i^\vee).
\eqno(14)
$$

As a basis $(e_1,\dots,e_n)$, consider now the image of $(\EE,\FF,\GG)$ in
K$_0(S)$. It is semiorthogonal since the corresponding collection of sheaves
is exceptional.

\smallskip
{\bf 3.2. Proposition. }{\it The image in\/ {\rm K}$_0(S)$ of the collection
of objects of\/ $D^b(S)$,
$$
\sigma=(G_\gamma\otimes K[2],\dots,G_1\otimes K[2],
L^D_\EE F_\beta,\dots,L^D_\EE F_1,E_\alpha,\dots,E_1),
$$
is a basis dual to\/} $(e_1,\dots,e_n)$.

\smallskip
{\sc Proof}. The collection $\sigma$ is exceptional since it is obtained from
the initial collection $(\EE,\FF,\GG)$ by mutations in $D^b(S)$. Namely,
if one introduces uniform indexing for the sheaves of the initial collection,
i.e.,
$$
(E_1,\dots,E_\alpha,F_1,\dots,F_\beta,G_1,\dots,G_\gamma)=
(A_1,\dots,A_n)\subset D^b(S),
$$
then the completeness, helix properties (Sec.~1.11), and triviality of
intrablock mutations imply that
$$
\sigma=(L_{A_1}\dots L_{A_{n-1}}A_n,\dots,L_{A_1}L_{A_{2}}A_{3},
L_{A_1}A_2,A_1)\triangleq(A_n^\vee,\dots,A_1^\vee).
$$
Since $\sigma$ is exceptional, its image in K$_0(S)$ is semiorthogonal.

It remains to show that $\chi(A_i,A_j^\vee)=\delta_{ij}$.

Consider the collections of objects of $D^b(S)$,
$$
(G_\gamma\otimes K[2],\dots,G_1\otimes K[2],
E_1,\dots,E_\alpha,F_1,\dots,F_\beta)
$$
and
$$
(L^D_\EE F_\beta,\dots,L^D_\EE F_1,E_1,\dots,E_\alpha,
G_1,\dots,G_\gamma).
$$
One easily sees that they are obtained by mutations from the initial
collection $(\EE,\FF,\GG)$ and, therefore, are exceptional. Hence,
$$
\begin{array}{l@{\qquad}l@{\qquad}l}
\chi(E_i,L^D_\EE F_j)=0, & \chi(F_i,E_j)=0, & \chi(G_i,L^D_\EE F_j)=0, \\
\chi(E_i,G_j\otimes K[2])=0, & \chi(F_i,G_j\otimes K[2])=0, &
\chi(G_i,E_j)=0, \\
\chi(E_i,E_j)=\delta_{ij}. &&
\end{array}
$$
The latter equality easily follows from the fact that $\EE$ is a block. By
Corollary 2.5 and equalities (11), $\chi(F_i,L^D_\EE F_j)=\delta_{ij}$. By
the Serre duality, $\chi(G_i,G_j\otimes K[2])=\chi(G_j,G_i)=\delta_{ij}$,
which completes the proof.

\medskip
{\bf 3.3. Derivation of Markov-type equations. }Introduce the notations
\begin{eqnarray*}&
x=r(\EE),\qquad y=r(\FF),\qquad z=r(\GG),\qquad y'=r(L^D_\EE\FF); &\\&
a=\chi(\FF,\GG),\qquad b=\chi(\GG\otimes K[2],\EE)=
\chi(\GG\otimes K,\EE),\qquad c=\chi(\EE,\FF).&
\end{eqnarray*}
 From the exact sequences of Proposition 2.2 and from Corollary 2.5,
$$
y'=y-c\alpha x,\qquad d(L^D_\EE\FF)=d(\FF)-c\alpha d(\EE).
$$

Let us compute values of the bilinear form (13) with the help of the
representation (14). As a basis in K$_0(S)$, let us use the image of the
initial collection $(\EE,\FF,\GG)$; and as a dual basis, the image of
$\sigma$ from the preceding proposition.

The equality $\ll r,r\rr=0$ means that $\alpha x^2+\beta yy'+\gamma z^2=0$,
or
$$
\alpha x^2+\beta y^2+\gamma z^2=c\alpha\beta xy.
\eqno(15)
$$
By the assumption, $x,y,z\ge0$ and $\alpha,\beta,\gamma>0$. The ranks $x$,
$y$, and $z$ cannot be zero simultaneously since the image of
$(\EE,\FF,\GG)$ is a basis in K$_0(S)$. Hence, $c>0$ and $x,y,z>0$.

Next, rewrite the equality $\ll r,d\rr-\ll d,r\rr=0$ as follows:
\begin{eqnarray*}
\lefteqn{\alpha xd(\EE)+\beta yd(L^D_\EE\FF)+
\gamma zd(\GG\otimes K[2])-\alpha d(\EE)x-\beta d(\FF)y'-\gamma zd(\GG)\;=}\\
 & = & \beta\Big(y(d(\FF)-c\alpha d(\EE))-d(\FF)(y-c\alpha x)\Big)+
\gamma z\Big(d(\GG\otimes K)-d(\GG)\Big)\;= \\
 & = & c\alpha\beta(d(\FF)x-d(\EE)y)-\gamma z^2K^2\;=\;
c^2\alpha\beta-\gamma z^2K^2\;=\;0.
\end{eqnarray*}
Hence,
$$
c=z\sqrt{\frac{K^2\gamma}{\alpha\beta}}.
\eqno(16)
$$

Substituting this into (15), we arrive at the statement below.

\smallskip
{\bf Theorem. }{\it The ranks $x,y,z$ and the numbers $\alpha,\beta,\gamma$
of sheaves in blocks of a complete collection\/ $(\EE,\FF,\GG)$ over a Del
Pezzo surface satisfy the relation
$$
\alpha x^2+\beta y^2+\gamma z^2=\sqrt{K^2\alpha\beta\gamma}\,xyz,
\eqno(3)
$$
where $K^2$ is a square of the canonical class of the surface, and the
coefficient in the right-hand side is an integer.}

\smallskip
Let us explain the latter assertion. As we have seen above, $x,y,z>0$, whence
$\sqrt{K^2\alpha\beta\gamma}\in{\Bbb Q}$. All rooted factors are integer;
hence, $\sqrt{K^2\alpha\beta\gamma}\in{\Bbb Z}$.

Furthermore, applying (16) to 3-block sheaf foundations
$(\FF,\GG,\EE\otimes(-K))$ and $(\GG\otimes K,\EE,\FF)$, we obtain
$$
a=x\sqrt{\frac{K^2\alpha}{\beta\gamma}}\qquad\hbox{and}\qquad
b=y\sqrt{\frac{K^2\beta}{\alpha\gamma}}.
\eqno(17)
$$
Then, expressing $x$, $y$, and $z$ in terms of $a$, $b$, and $c$,
and substituting these expressions into (3), we get the equation on the
dimensions of the interblock \Hom\ spaces as follows:
$$
\frac{a^2}\alpha+\frac{b^2}\beta+\frac{c^2}\gamma=abc.
$$

\medskip
{\bf 3.4. Corollary. }1. {\it For any pair of blocks $(\EE,\FF)$ contained in
a\/ $3$-block collection,\/ $\Hom(E_i,F_j)\ne0$. In other words, there is no
complete two-block collection over a Del Pezzo surface.

$2$. Any mutation of a complete\/ $3$-block collection is a division.

$3$. Any sheaf contained in a complete\/ $3$-block collection is locally
free.}

\smallskip
{\sc Proof}. Item 1 follows from the fact that
$c=\chi(\EE,\FF)>0$. Then, in a pair $(\EE,\FF)$, a mutation of the extension
type can never occur. Hence, an inverse mutation of the recoil type is also
impossible. Item 3 follows from the classification of exceptional sheaves
(Sec.\ 1.3) and the inequalities $x,y,z>0$ obtained above.

\medskip
{\bf 3.5. The list of equations. }It is well known [14], that a Del Pezzo
surface $S$ is isomorphic to either $\PP^1\times\PP^1$ or $X_m$, the latter
being a plane with $m$ generic points blown-up, where $0\le m\le8$.

Assume there is a complete 3-block collection $(\EE,\FF,\GG)$ over $S$. Then,
as above, the ranks $x$, $y$, and $z$ are positive and obey the equation (3).
Moreover, the image of a given collection in K$_0(S)$ forms a basis, whence
$\alpha+\beta+\gamma={\rm rk\,K}_0(S)=12-K^2$. Therefore, $\alpha$,
$\beta$, $\gamma$, and the squared anticanonical class of $S$ satisfy the
following system of conditions:
$$
\left\{\begin{array}{l}
\alpha+\beta+\gamma+K^2=12,\\
K^2\alpha\beta\gamma\hbox{ is a square of an integer},\\
\alpha,\beta,\gamma\ge1,\\
1\le K^2\le9.
\end{array}\right.
$$

The system is solved by the finite exhausting. As an answer, we present the
complete list of equations (3) with indication of the correspondent
surfaces\footnote{We prove below that for each equation, the corresponding
complete 3-block collection of sheaves exists.}. For some equations, a common
multiplier can be cancelled, but we do not do this in order to leave as
coefficients the numbers $\alpha$, $\beta$, and $\gamma$ of sheaves in
the correspondent blocks. One can easily check that the blocks with given
numbers of sheaves can be arbitrarily reordered by mutations. (For instance,
the first block of $(L_\EE\FF,\EE,\GG)$ consists of $\beta$ sheaves, the
second one consists of $\alpha$ sheaves, and the third one, of $\gamma$
sheaves\@.) Therefore, without loss of generality, we assume that
$\alpha\le\beta\le\gamma$.

\smallskip
{\bf Definition. }A {\it minimum solution\/} of the equation (3) is a
positive integer solution with the minimum possible sum $x+y+z$.

\smallskip
{\bf Proposition. }{\it All equations of type\/ $(3)$ with
$\alpha\le\beta\le\gamma,$ together with their minimum solutions, are
presented in the table below.}
%\vskip \abovedisplayskip
\begin{center}\renewcommand{\arraystretch}{1.4} 
\begin{tabular}{|c|c|c|l|}\hline 
Number of the & Surface & Equation              & Minimum     \\
equation  &             &                       & solution    \\ \hline
\x1       & $\PP^2$     & $ x^2+ y^2+ z^2=3xyz$ & $(1,1,1)$   \\
\x2       & $\PP^1
           \times\PP^1$ & $ x^2+ y^2+2z^2=4xyz$ & $(1,1,1)$   \\
\x3       & $X_3$       & $ x^2+2y^2+3z^2=6xyz$ & $(1,1,1)$   \\
\x4       & $X_4$       & $ x^2+ y^2+5z^2=5xyz$ & $(1,2,1)$ and
                                                  $(2,1,1)$   \\
\x5       & $X_5$       & $2x^2+2y^2+4z^2=8xyz$ & $(1,1,1)$   \\
\x{6.1}   & $X_6$       & $3x^2+3y^2+3z^2=9xyz$ & $(1,1,1)$   \\
\x{6.2}   & $X_6$       & $ x^2+2y^2+6z^2=6xyz$ & $(2,1,1)$   \\
\x{7.1}   & $X_7$       & $ x^2+ y^2+8z^2=4xyz$ & $(2,2,1)$   \\
\x{7.2}   & $X_7$       & $2x^2+4y^2+4z^2=8xyz$ & $(2,1,1)$   \\
\x{7.3}   & $X_7$       & $ x^2+3y^2+6z^2=6xyz$ & $(3,1,1)$   \\
\x{8.1}   & $X_8$       & $ x^2+ y^2+9z^2=3xyz$ & $(3,3,1)$   \\
\x{8.2}   & $X_8$       & $ x^2+2y^2+8z^2=4xyz$ & $(4,2,1)$   \\
\x{8.3}   & $X_8$       & $2x^2+3y^2+6z^2=6xyz$ & $(3,2,1)$   \\
\x{8.4}   & $X_8$       & $ x^2+5y^2+5z^2=5xyz$ & $(5,2,1)$ and
                                                  $(5,1,2)$   \\ \hline
\end{tabular}
\end{center}
%\vskip \abovedisplayskip

\smallskip
{\sc Proof}. We omit verification of the fact that all possible equations of
type (3) are presented here. Let us show that the right-hand column actually
contains minimum solutions.

For the equations that have a solution $(1,1,1)$, this is obvious. For the
equations \x4, \x{6.2}, and \x{7.2}, this is true since $(1,1,1)$ is not a
solution to these equations. Consider the remaining cases.

\x{7.1}: It is easily seen that $x$ and $y$ are of the same parity, but they
cannot be odd since $x^2+y^2$ is divisible by 4.

\x{7.3}: $x$ is divisible by 3.

\x{8.1}: $x^2+y^2$ is divisible by 3, but since a square's residual
modulo 3 equals either 0 or 1, both $x$ and $y$ are divisible by 3.

\x{8.2}: $x$ is even; hence, $2y^2$ is divided by 4. Then $y$ is even, and
hence, $x^2$ is divisible by 8, i.e., $x$ is also divisible by 4.

\x{8.3}: $x$ is divisible by 3, and $y$ is even.

\x{8.4}: $x$ is divisible by 5. Let $x=5\widetilde x$, then the equation
takes the form $5\widetilde x^2+y^2+z^2=5\widetilde xyz$, which coincides with
\x4 up to a designation.

\medskip
{\it Remark}. The surface $X_1$ has $K^2=8$, but there are no complete
exceptional collections over $X_1$. Indeed, the total number of sheaves in a
complete collection should be equal to ${\rm rk\,K}_0(X_1)=4$, i.e., one
block should consist of two sheaves. By Proposition 1.6, the difference
$c$ of the first Chern classes of these sheaves should satisfy the relations
$c^2=-2$ and $c\cdot K=0$. But there is no such a divisor over $X_1$.

\medskip
{\bf 3.6. Solution mutations. }Let $(x,y,z)$ be a solution of (3).
The {\it solution mutation in the variable\/} $y$ is the mapping
$M_y\colon (x,y,z)\mapsto (x,y',z)$, where
$$
y'=\sqrt{\frac{K^2\alpha\gamma}\beta}xz-y=c\alpha x-y=a\gamma z-y.
$$
Similarly, the solution mutations $M_x$ and $M_z$ (in $x$ and $z$
respectively) are defined.

For any complete 3-block collection with structure $\{\alpha,
\beta,\gamma\}$, define the {\it correspondent solution\/} of (3) as a
triple of numbers $(r_\alpha,r_\beta,r_\gamma)$, where the first number is
the rank of sheaves in the block of length $\alpha$, etc. If $\alpha < \beta
<\gamma$, a single solution corresponds to the 3-block collection, and for
$\alpha=\beta$ or $\beta =\gamma$, more than one solution may correspond to
one collection. Moreover, the same solution may correspond to collections of
different types obtained by various permutations of $(\alpha, \beta,\gamma)$.

Let $(x,y,z)$ be a solution corresponding to a collection $(\EE,\FF,\GG)$ of
type $(\alpha, \beta,\gamma)$. Then the solution $(x,y',z)$ described above
corresponds to the collections $(L_\EE\FF,\EE,\GG)$ and $(\EE,\GG,R_\GG\FF)$.
This follows from Corollary 3.4 and Proposition 2.2. Note that
$L_\EE\FF=R_\GG\FF\otimes K$. In the general case, one can easily check that
mutations of the collection which change the block of length
$\alpha$ ($\beta$ or $\gamma$) induce the correspondent solution mutations
in $x$ ($y$ or $z$).

Thus, mutations of complete 3-block collections are concordant with mutations
of correspondent solutions.

\medskip
{\bf 3.7. Proposition. }(a) {\it Any solution of any of equation from
Proposition\/ $3.5$ can be reduced by mutations to a minimum solution.}

(b) {\it Moreover, for a nonminimum solution $(x,y,z)$, one mutation reduces
and two others increase the sum\/} $x+y+z$.

\smallskip
{\sc Proof}. For the equations \x1 and \x2, (a) is well known; see [15, 18,
19]. Let us verify (b) for these equations. Since \x1 is symmetric in $x$,
$y$, and $z$, we may assume that $x=\max\{x,y,z\}$. Then $y+y'=3xz\ge 3yz\ge
3y>2y$, whence $y'>y$, i.e., the mutation in $y$ increases $x+y+z$. For the
mutation in $z$, the reasoning is similar.

The equation \x2 is symmetric in $x,y$. Put, for definiteness,
$x\ge y$. Consider the case $x\ge z$. Then
$y+y'=4xz\ge 4yz\ge 4y>2y\ \Longrightarrow\ y'>y$, and also
$z+z'=2xy\ge 2zy\ge 2z\ \Longrightarrow\ z'\ge z$.
The equality $z'=z$ is possible only if $x=z$ and $y=1$. Putting this into
\x2, we get $x=1$. Hence, if $(x,y,z)\ne (1,1,1)$, then $z'>z$. For the case
$x<z$, the reasoning is similar.

Next, consider the equations \x3 and \x4. Introduce the notations
$$
\Phi(x,y,z)=\alpha x^2+\beta y^2+\gamma
z^2-\sqrt{K^2\alpha\beta\gamma}\,xyz,\qquad\varphi_y(t)=\Phi(x,t,z),
$$
$(x,y,z)$ being a fixed solution of (3). Here, $t_1=y$ and
$t_2=y'$ are the roots of the quadratic equation $\varphi_y(t)=0$. Consider
also the functions $\varphi_x(t)=\Phi(t,y,z)$ and $\varphi_z(t)=\Phi(x,y,t)$.
The main tool for proving (a) is the following obvious statement:
$$
y'\ge y\quad\Longleftrightarrow\quad \varphi_y(t)\ge0,\  \forall t\le y,
$$
and the analogous statements for $x$ and $z$ as well.

The condition $y'\ge y$ means that the mutation of the solution $(x,y,z)$
in $y$ does not decrease $y$.

Let $(x,y,z)$ be a solution of \x3 such that none of its mutation reduces
$x+y+z$. Let us show that $(x,y,z)=(1,1,1)$ then. Consider the cases below.

(i) $x\ge y\ge z$. Then $0\le\varphi_x(y)=3y^2+3z^2-6y^2z\le6y^2-6y^2z$,
whence $z=1$ and $\varphi_x(y)=3-3y^2\ge0$. Therefore, $y=1$, and \x3
directly implies $x=1$.

(ii) $x\ge z>y$. Then $0\le\varphi_x(z)=4z^2+2y^2-6z^2y<6z^2-6z^2y$, which is
impossible. The other cases are also impossible:

(iii)
$y>x\ge z\;\Longrightarrow\;0<\varphi_y(x)=3x^2+3z^2-6x^2z\le6x^2-6x^2z$;

(iv) $z>x\ge y\;\Longrightarrow\;0<\varphi_z(x)=4x^2+2y^2-6x^2y\le6x^2-6x^2y$;

(v) $y\ge z>x\;\Longrightarrow\;0<\varphi_y(z)=x^2+5z^2-6z^2x<6z^2-6z^2x$;

(vi) $z>y>x\;\Longrightarrow\;0<\varphi_z(y)=x^2+5y^2-6y^2x<6y^2-6y^2x$.

Thus, for any nonminimum solution of \x3, a mutation exists which reduces
$x+y+z$. Hence, for \x3, (a) is proved. Let us prove (b) under
$x=\max\{x,y,z\}$ (in other cases, the reasoning is similar). We have
$y+y'=3xz\ge 3yz\ge 3y>2y\ \Longrightarrow\ y'>y$, and also
$z+z'=2xy\ge 2zy\ge 2z\ \Longrightarrow\ z'\ge z$.
The equality $z'=z$ is possible only if $x=z$ and $y=1$, whence the
minimality of $(x,y,z)$ follows.

Now, let $(x,y,z)$ be a solution of \x4 such that none of its mutations
reduces $x+y+z$. Let us show that $(x,y,z)$ coincides with one of the minimum
solutions, $(2,1,1)$ or $(1,2,1)$. The variables $x$ and $y$ are equivalent,
so we assume $x\ge y$.

(i) $x\ge y\ge z$. Then $0\le\varphi_x(y)=2y^2+5z^2-5y^2z\le7y^2-5y^2z$,
whence $z=1$ and $\varphi_x(y)=5-3y^2\ge0$. Hence, $y=1$, and then $x=2$,
which directly follows from \x4.

(ii) $x\ge z>y$. Then $0\le\varphi_x(z)=6z^2+y^2-5z^2y<7z^2-5z^2y$, whence
$y=1$. Consider the solution mutation in $z$, $z'=xy-z=x-z$. By the
assumption, $z'\ge z$; hence, $x\ge2z$. Then
$0\le\varphi_x(2z)=9z^2+1-10z^2=1-z^2$, i.e., $z=1$, which provides a
contradiction.

(iii) $z>x\ge y$. Then $0\le\varphi_z(x)=6x^2+y^2-5x^2y\le7x^2-5x^2y$, whence
$y=1$, and the solution mutation in $z$ yields $z'=x-z<0$, a contradiction.

Thus, for any nonminimum solution of \x4, a mutation exists which reduces
$x+y+z$. This proves (a). Let us prove (b), assuming $x\ge y$ as before.

(i) $x\ge z$. Then
$y+y'=5xz\ge 5yz\ge 5y>2y\ \Longrightarrow\ y'>y$, i.e., the mutation
in $y$ increases $x+y+z$. Let $y\ge 2$, then $z+z'=xy\ge 2x\ge 2z$, whence
$z'\ge z$. Here, the equality is possible only if
$y=2$ and $x=z$, which implies $x=z=1$. Hence, the mutation in $z$ reduces
$z$ if $(x,y,z)\ne (1,2,1)$. Now, let $y=1$. Let us show that the mutations
in $x$ and $z$ cannot reduce $x+y+z$ simultaneously. Indeed, if $x'=5z-x<x$,
then $5z<2x$, whence $2z<x$, and $z'=x-z<z$.

(ii) $z>x$. In this case, $y+y'=5xz\ge 5xy\ge 5y>2y$, whence $y'>y$, i.e.,
the mutation in $y$ increases $x+y+z$. This completes the proof of
(b) for the equation \x4 due to its symmetry in $x,y$.

Thus, the proposition is valid for the first four equations. Each of the
others can be reduced to one of the first four by a change of variables,
which is possible in each case (according to the proof of Proposition 3.5).
Namely,

\x{6.2} is reduced to \x3 under $x=2\widetilde x$;

\x{7.1} is reduced to \x2 under $x=2\widetilde x$, $y=2\widetilde y$;

\x{7.2} is reduced to \x2 under $x=2\widetilde x$;

\x{7.3} is reduced to \x3 under $x=3\widetilde x$;

\x{8.1} is reduced to \x1 under $x=3\widetilde x$, $y=3\widetilde y$;

\x{8.2} is reduced to \x2 under $x=4\widetilde x$, $y=2\widetilde y$;

\x{8.3} is reduced to \x3 under $x=3\widetilde x$, $y=2\widetilde y$;

\x{8.4} is reduced to \x4 under $x=5\widetilde x$.

This completes the proof of the proposition.

\medskip
{\bf 3.8. Groups of equations. } As one can see from the proof of the
preceding proposition, any equation starting from \x5 either can be obtained
from one of the first four equations by a change of variables, or is
proportional to one of them. Let us join the equations in groups as follows:
\begin{center}
\begin{tabular}{rl@{$\qquad$}rl}
Group I:& \x1, \x{6.1}, \x{8.1};&
Group II:& \x2, \x5, \x{7.1}, \x{7.2}, \x{8.2};\\
Group III:& \x3, \x{6.2}, \x{7.3}, \x{8.3};&
Group IV:& \x4, \x{8.4}.
\end{tabular}
\end{center}
For each equation, consider a pseudograph whose vertices are the solutions,
and two vertices are joined by an edge if and only if the solutions can be
obtained one from another by one of the mutations $M_x$, $M_y$, or $M_z$.
Obviously, the pseudographs for equations of the same group are isomorphic.
One can easily deduce from the preceding proposition that the pseudographs of
the solutions have the form as follows: for the group I, $\Gamma_1$ (see
Fig.~1);\setlength{\unitlength}{0.01in}%
\sbox{\hdts}{%                        %
\begin{picture}(14,3)                 %
\multiput(1,1)(6,0){3}{\circle*{2.5}} %
\end{picture}%                        %
}%
\sbox{\ldts}{%                        %
\begin{picture}(8,12)                 %
\multiput(-3,-5.19615)(3,5.19615){3}% %
{\circle*{2.5}}                       %
\end{picture}%
}%
\sbox{\rdts}{%                        %
\begin{picture}(8,12)                 %
\multiput(-3,5.19615)(3,-5.19615){3}% %
{\circle*{2.5}}                       %
\end{picture}%
}%
\begin{figure}
\begin{center}
\begin{picture}(520,220)
\put(108,0){\Large $\Gamma_1$}
\put(403,40){\Large $\Gamma_2$}
\put(0,30){%
\begin{picture}(225,185)
\special{em:linewidth 0.008in}
\put(104.5833,114){$P_0$}
\put(112.5833,104){\circle*{4}}  %
\put(112.5833,50){\circle*{4}}   %
\put(65.8179,131){\circle*{4}}   %
\put(159.3486,131){\circle*{4}}  %
\put(65.8179,167){\circle*{4}}   %
\put(159.3486,167){\circle*{4}}  %
\put(34.6410,113){\circle*{4}}   %
\put(190.5256,113){\circle*{4}}  %
\put(81.4064,32){\circle*{4}}    %
\put(143.7602,32){\circle*{4}}   %
\put(112.5833,104){\special{em:moveto}}  %
\put(112.5833,50){\special{em:lineto}}   %
\put(65.8179,131){\special{em:moveto}}   %
\put(112.5833,104){\special{em:lineto}}  %
\put(159.3486,131){\special{em:moveto}}  %
\put(112.5833,104){\special{em:lineto}}  %
\put(65.8179,167){\special{em:moveto}}   %
\put(65.8179,131){\special{em:lineto}}   %
\put(159.3486,167){\special{em:moveto}}   %
\put(159.3486,131){\special{em:lineto}}   %
\put(65.8179,131){\special{em:moveto}}   %
\put(34.6410,113){\special{em:lineto}}   %
\put(159.3486,131){\special{em:moveto}}   %
\put(190.5256,113){\special{em:lineto}}   %
\put(112.5833,50){\special{em:moveto}}   %
\put(81.4064,32){\special{em:lineto}}   %
\put(112.5833,50){\special{em:moveto}}   %
\put(143.7602,32){\special{em:lineto}}   %
\put(47.1769,179){\special{em:moveto}}   %
\put(65.8179,167){\special{em:lineto}}   %
\put(88.7461,179){\special{em:moveto}}   %
\put(65.8179,167){\special{em:lineto}}   %
\put(41.1143,182.5){\usebox{\ldts}}
\put(94.8087,182.5){\usebox{\rdts}}
\put(138.5640,179){\special{em:moveto}}   %
\put(159.3486,167){\special{em:lineto}}   %
\put(180.1333,179){\special{em:moveto}}   %
\put(159.3486,167){\special{em:lineto}}   %
\put(132.5014,182.5){\usebox{\ldts}}
\put(186.1959,182.5){\usebox{\rdts}}
\put(81.4064,32){\special{em:moveto}}   %
\put(81.4064,8){\special{em:lineto}}   %
\put(60.6218,44){\special{em:moveto}}   %
\put(81.4064,32){\special{em:lineto}}   %
\put(74.4064,0){\usebox{\hdts}}
\put(54.5592,47.5){\usebox{\ldts}}
\put(143.7602,32){\special{em:moveto}}   %
\put(143.7602,8){\special{em:lineto}}   %
\put(164.5445,44){\special{em:moveto}}   %
\put(143.7602,32){\special{em:lineto}}   %
\put(136.7602,0){\usebox{\hdts}}
\put(170.6071,47.5){\usebox{\rdts}}
\put(211.3102,125){\special{em:moveto}}   %
\put(190.5256,113){\special{em:lineto}}   %
\put(190.5256,113){\special{em:moveto}}   %
\put(190.5256,89){\special{em:lineto}}   %
\put(183.5256,81){\usebox{\hdts}}
\put(217.3831,128.5){\usebox{\rdts}}
\put(13.8564,125){\special{em:moveto}}   %
\put(34.6410,113){\special{em:lineto}}   %
\put(34.6410,113){\special{em:moveto}}   %
\put(34.6410,89){\special{em:lineto}}   %
\put(27.6410,81){\usebox{\hdts}}
\put(7.7938,128.5){\usebox{\ldts}}
\end{picture}
}
\put(295,30){%
\begin{picture}(225,185)
\special{em:linewidth 0.008in}
\put(104.5833,114){$P_0$}
\put(112.5833,104){\circle*{4}}  %
\put(65.8179,131){\circle*{4}}   %
\put(159.3486,131){\circle*{4}}  %
\put(65.8179,167){\circle*{4}}   %
\put(159.3486,167){\circle*{4}}  %
\put(34.6410,113){\circle*{4}}   %
\put(190.5256,113){\circle*{4}}  %
\put(65.8179,131){\special{em:moveto}}   %
\put(112.5833,104){\special{em:lineto}}  %
\put(159.3486,131){\special{em:moveto}}  %
\put(112.5833,104){\special{em:lineto}}  %
\put(65.8179,167){\special{em:moveto}}   %
\put(65.8179,131){\special{em:lineto}}   %
\put(159.3486,167){\special{em:moveto}}   %
\put(159.3486,131){\special{em:lineto}}   %
\put(65.8179,131){\special{em:moveto}}   %
\put(34.6410,113){\special{em:lineto}}   %
\put(159.3486,131){\special{em:moveto}}   %
\put(190.5256,113){\special{em:lineto}}   %
\put(47.1769,179){\special{em:moveto}}   %
\put(65.8179,167){\special{em:lineto}}   %
\put(88.7461,179){\special{em:moveto}}   %
\put(65.8179,167){\special{em:lineto}}   %
\put(41.1143,182.5){\usebox{\ldts}}
\put(94.8087,182.5){\usebox{\rdts}}
\put(138.5640,179){\special{em:moveto}}   %
\put(159.3486,167){\special{em:lineto}}   %
\put(180.1333,179){\special{em:moveto}}   %
\put(159.3486,167){\special{em:lineto}}   %
\put(132.5014,182.5){\usebox{\ldts}}
\put(186.1959,182.5){\usebox{\rdts}}
\put(211.3102,125){\special{em:moveto}}   %
\put(190.5256,113){\special{em:lineto}}   %
\put(190.5256,113){\special{em:moveto}}   %
\put(190.5256,89){\special{em:lineto}}   %
\put(183.5256,81){\usebox{\hdts}}
\put(217.3831,128.5){\usebox{\rdts}}
\put(13.8564,125){\special{em:moveto}}   %
\put(34.6410,113){\special{em:lineto}}   %
\put(34.6410,113){\special{em:moveto}}   %
\put(34.6410,89){\special{em:lineto}}   %
\put(27.6410,81){\usebox{\hdts}}
\put(7.7938,128.5){\usebox{\ldts}}
\put(102.5233,56.52)%                    %
{%                                       %
\begin{picture}(20.13,47.48)             %
\special{em:linewidth 0.008in}
\put(26,20){$M$}                         %
\put(10.06,47.48){\special{em:moveto}}
\put(8.10,43.75){\special{em:lineto}}
\put(6.36,40.18){\special{em:lineto}}
\put(4.86,36.76){\special{em:lineto}}
\put(3.60,33.50){\special{em:lineto}}
\put(2.53,30.38){\special{em:lineto}}
\put(1.66,27.41){\special{em:lineto}}
\put(0.98,24.60){\special{em:lineto}}
\put(0.50,21.95){\special{em:lineto}}
\put(0.18,19.43){\special{em:lineto}}
\put(0.01,17.10){\special{em:lineto}}
\put(0.00,14.90){\special{em:lineto}}
\put(0.13,12.85){\special{em:lineto}}
\put(0.40,10.95){\special{em:lineto}}
\put(0.78,9.20){\special{em:lineto}}
\put(1.26,7.60){\special{em:lineto}}
\put(1.86,6.16){\special{em:lineto}}
\put(2.55,4.86){\special{em:lineto}}
\put(3.31,3.73){\special{em:lineto}}
\put(4.15,2.73){\special{em:lineto}}
\put(5.05,1.90){\special{em:lineto}}
\put(5.98,1.21){\special{em:lineto}}
\put(6.98,0.68){\special{em:lineto}}
\put(8.00,0.30){\special{em:lineto}}
\put(9.03,0.08){\special{em:lineto}}
\put(10.06,0.00){\special{em:lineto}}
\put(11.11,0.08){\special{em:lineto}}
\put(12.15,0.30){\special{em:lineto}}
\put(13.16,0.68){\special{em:lineto}}
\put(14.15,1.21){\special{em:lineto}}
\put(15.08,1.90){\special{em:lineto}}
\put(15.98,2.73){\special{em:lineto}}
\put(16.81,3.73){\special{em:lineto}}
\put(17.58,4.86){\special{em:lineto}}
\put(18.26,6.16){\special{em:lineto}}
\put(18.86,7.60){\special{em:lineto}}
\put(19.35,9.20){\special{em:lineto}}
\put(19.73,10.95){\special{em:lineto}}
\put(20.00,12.85){\special{em:lineto}}
\put(20.13,14.90){\special{em:lineto}}
\put(20.11,17.10){\special{em:lineto}}
\put(19.95,19.43){\special{em:lineto}}
\put(19.63,21.95){\special{em:lineto}}
\put(19.15,24.60){\special{em:lineto}}
\put(18.46,27.41){\special{em:lineto}}
\put(17.60,30.38){\special{em:lineto}}
\put(16.53,33.50){\special{em:lineto}}
\put(15.26,36.76){\special{em:lineto}}
\put(13.76,40.18){\special{em:lineto}}
\put(12.05,43.75){\special{em:lineto}}
\put(10.06,47.48){\special{em:lineto}}
\end{picture}%
}
\end{picture}
}
\end{picture}\\
Fig.\ 1.
\end{center}
\end{figure}
for the groups II and III, $\Gamma_2$; and for the group IV, the pseudograph
consists of two connected components isomorphic to $\Gamma_2$. The point
$P_0$ denotes the minimum solution. Moreover, $\Gamma_1$ and $\Gamma_2$ have
no cycles, i.e., $\Gamma_1$ is actually a graph. By abuse of language, we
will say ``solution graph'' instead of ``solution pseudograph'' although
$\Gamma_2$ is not a graph (since it has a single loop $M$ starting and ending
at the minimum solution).

\medskip
{\bf 3.9. Definition. }A complete 3-block collection $(\EE,\FF,\GG)$ is
called {\it minimum\/} if the sum $r(\EE)+r(\FF)+r(\GG)$ is minimum of all
sums for 3-block collections of the same structure.

\smallskip
{\bf 3.10. Theorem. }{\it Any complete\/ $3$-block collection can be obtained
by mutations from a minimum collection of type $(\alpha,\beta,\gamma)$, where
$\alpha\le\beta\le\gamma$.}

\smallskip
{\sc Proof}. According to 3.6 and 3.7, any complete 3-block collection can be
reduced by mutations to a minimum one. Invertibility of mutations imply that
a given collection can be obtained by mutations from a minimum one. Hence, it
suffices to show that a minimum collection whose type is an arbitrary
permutation of $(\alpha,\beta,\gamma)$ can be obtained from a minimum
collection $(\EE,\FF,\GG)$ of type $(\alpha,\beta,\gamma)$, where
$\alpha\le\beta\le\gamma$. The sequences of mutations $L_1L_2$ and $R_2R_1$
(cf.\ the notations of Sec.~2.6) take $(\EE,\FF,\GG)$ to $(\GG(K),\EE,\FF)$
and $(\FF,\GG,\EE(-K))$ respectively. The types of these collections are
cyclic transpositions of $(\alpha,\beta,\gamma)$. Thus, the theorem is valid
for $\alpha=\beta$ or $\beta=\gamma$. This condition holds for all equations
except those of group III. For the equations of this group, the mutation that
preserves the minimum solution is induced by the mutation of the minimum
collection which performs a transposition of $(\alpha,\beta,\gamma)$. The
theorem is proved.

\smallskip
Thus, the problem on the action of the braid group and, in particular, on the
set of orbits, reduces to the problem on the set of minimum collections.

Below, we study minimum 3-block collections, namely, prove their existence
and describe the action of the Weyl group on collections of a given
structure.

\bigskip
\begin{center}{\large\bf 4. Existence of minimum collections.}\end{center}

\medskip
{\bf 4.1. Agreements. }Fix generic points $x_1,\dots,x_8$ on a projective
plane and denote by $\sigma_r\colon X_r\rightarrow\PP^2$ the monoidal
transform with center $\{x_1,\dots,x_r\}$. Let $\ell_0$ be a divisor class on
$X_r$ equal to the lifting of the class of a line on $\PP^2$. Let $\ell_i$,
$i=1,\dots,r$, be the classes of exceptional curves $\sigma_r^{-1}(x_i)$. For
$1\le i\le r$, each of $\ell_i$'s contains a single divisor, namely, the curve
$\sigma_r^{-1}(x_i)$ itself. We denote this curve by the same symbol $\ell_i$.

The surface $X_{r+1}$ is obtained as $X_r$ with a blown-up point
$\sigma_r^{-1}(x_{r+1})$, so we have the diagram
$$
X_8\longrightarrow X_7\longrightarrow\dots\longrightarrow
  X_2\longrightarrow X_1\longrightarrow X_0=\PP^2.
\eqno(18)
$$
In correspondence with it, we will consider $\ZZ$-modules $\pic X_r$ to be
embedded in one another, i.e.,
$$
\ZZ\cong\pic\PP^2\subset\pic X_1\subset\dots\subset\pic X_8\ ,
$$
where $\pic X_r=\ZZ \ell\oplus\ZZ \ell_1\oplus\dots\oplus\ZZ \ell_r$. The intersection
form is defined by
$$
\ell_0^2=1 ,\qquad \ell_i^2=-1,\quad i\ge 1,\qquad \ell_i\cdot \ell_j=0,\quad i\ne j\ .
$$
The canonical class of $X_r$ is
$$
\omega_r=-3\ell+\sum_{i=1}^r \ell_i .
$$

Denote by $\sigma\colon X_r\rightarrow X_p$, $r>p$, the composition
of morphisms in (18), i.e., blowing up the points
$\sigma_p^{-1}(x_{p+1}),\dots ,\sigma_p^{-1}(x_r)$. This will not cause
ambiguity since it will always be clear from a context which $p$ and
$r$ are meant. Under our agreements, $\pic X_p\subset\pic X_r$ is the
orthogonal complement to the linear span of $\ell_{p+1},\dots,\ell_r$.

To any divisor class $a\ell_0+\sum\limits_{i=1}^r b_i\ell_i$ modulo linear
equivalence, a unique class of invertible sheaves modulo isomorphism
corresponds which we denote by
$\co_{X_r}\Big(a\ell_0+\sum\limits_{i=1}^rb_i\ell_i\Big)$. Here,
$$
\sigma^{\ast}\co_{X_p}\Big(a\ell_0+\sum_{i=1}^pb_i\ell_i\Big)=
\co_{X_r}\Big(a\ell_0+\sum_{i=1}^pb_i\ell_i\Big).
$$

In this section, we show that over the surfaces $X_r,$ $3\le r\le 8$,
complete 3-block exceptional collections corresponding to minimum solutions
of Markov-type equations (see Sec.\ 3.5) exist, and all such collections can
be obtained by a procedure which may be called ``lifting.'' The
general scheme of this procedure is as follows. It is known that if $\tau$ is
a complete exceptional collection over $X_p$, then the collection
$\sigma^{\ast}\tau$ complemented from the left by the sheaves
$\co_{\ell_{p+1}}(-1),\dots ,\co_{\ell_r}(-1)$ (or from the right by
$\co_{\ell_{p+1}},\dots ,\co_{\ell_r}$), form a complete exceptional collection
over $X_r$. Here, if $\tau$ is 3-block, the latter collection over $X_r$ is
4-block. In some cases%
\footnote{Namely, if the difference between the slopes of
two blocks in $\tau$ equals the slope of the anticanonical class of
$X_r$.} 
it turns out to be possible to perform such mutations of the 4-block
collection $(\sigma^{\ast}\tau,\co_{\ell_{p+1}},\dots,\co_{\ell_r})$ (or
$(\co_{\ell_{p+1}}(-1),\dots,\co_{\ell_r}(-1),\sigma^{\ast}\tau)$) that in an
obtained 4-block collection, two neighboring blocks can be joined into one,
which actually gives a 3-block collection over $X_r$.

\smallskip
Recall that a minimum collection over a Del Pezzo surface is a complete
3-block collection whose block ranks form a minimum solution of the
correspondent Markov-type equation.

\smallskip
{\bf 4.2. Proposition. }{\it Over the surfaces $X_r,$ $3\le r\le 8,$ there
exist minimum collections of types\/ $(\alpha,\beta,\gamma)$, where\/}
$\alpha\le\beta\le\gamma$.

\smallskip
{\sc Proof}. We could just present the collections required, but it will be
essential for us that all of them can be obtained by the ``lifting''
described above. Therefore, we present the corresponding sequences of
mutations using the notations of Sec.\ 2.6. For known discrete invariants
$r(E_i)$ and $\chi(E_i,E_j)$, $i<j$, Propositions 2.2 and 2.3 make it
possible to determine the type of any mutation of a 4-block collection
$({\cal E}_1,{\cal E}_2,{\cal E}_3,{\cal E}_4)$ and compute these invariants
for the collection obtained by a mutation. Note that if a collection
$({\cal E}_1,{\cal E}_2,{\cal E}_3,{\cal E}_4)$ is complete, the sequences of
mutations
$$
R_1^{(3)}=R_3\circ R_2\circ R_1\qquad\mbox{and}\qquad
L_3^{(3)}=L_1\circ L_2\circ L_3
$$
take it to $({\cal E}_2,{\cal E}_3,{\cal E}_4,{\cal E}_1(-K))$ and
$({\cal E}_4(K),{\cal E}_1,{\cal E}_2,{\cal E}_3)$ respectively.

We enumerate the items of the proof in the same way as the Markov-type
equations in the table of Sec.\ 3.5. For brevity, we denote a block
$\co_{\ell_{p+1}},\dots ,\co_{\ell_r}$, $p<r$, consisting of torsion sheaves by
$\lll p,r$, or by $\LL p,{p+1}$ for $r=p+1$. We arrange the final 3-block
collection into a table
$\begin{array}{|c|c|c|}
\hline
{\cal E}&{\cal F}&{\cal G}\\
\hline
\end{array}\,$, where the block obtained as a result of the last operation
(joining two blocks into one) is divided by a dotted line into the parts of
which it is composed. We call this block {\it distinguished}. Introduce the
notation
$$
\tau_0=(\co_{\PP^2}(-1),\co_{\PP^2},\co_{\PP^2}(1))
$$
for a well-known foundation of a helix over $\PP^2$.

\smallskip
\x3. Consider the sequence of mutations
$$
\settowidth{\plw}{$\ \co_{X_3}(2\ell-\ell_1-\ell_2-\ell_3)\ $}
   \begin{array}{c}\medskip
   (\sigma^{\ast}_3\tau_0,\lll 1,3)=
    (\co_{X_3}(-\ell),\,\co_{X_3},\,\co_{X_3}(\ell),\,\lll 1,3)
       \stackrel{ R_1^{(3)}}{\longrightarrow}
\\ \bigskip
\longrightarrow (\co_{X_3},\,\co_{X_3}(\ell),\,\lll 1,3,\,
\co_{X_3}(2\ell-\ell_1-\ell_2-\ell_3))
  \stackrel{R_3}{\longrightarrow}\\
\longrightarrow \renewcommand{\arraystretch}{1.3}
\begin{array}{|c|@{}c@{}|c|}
\hline
      &\raisebox{-0.2em}[0.8em][0em]{$\co_{X_3}(\ell)$}&\co_{X_3}(2\ell-\ell_2-\ell_3)\\
\ \co_{X_3}&\pline{0.3em},{\plw}&\co_{X_3}(2\ell-\ell_1-\ell_3)\\
         &\raisebox{0.2em}[0.4em][0.4em]{$\co_{X_3}(2\ell-\ell_1-\ell_2-\ell_3)$}
&\co_{X_3}(2\ell-\ell_1-\ell_2)\\
\hline
\end{array}=\tau\xx3.
\end{array}
$$
This is the desired collection over a plane with three blown-up points.

\smallskip
\x4. Over $X_4$, we have
$$
\settowidth{\plw}{$\ \co_{X_4}(\ell_4-\omega_4-\ell)\ $}
  \begin{array}{c}\medskip
   (\sigma^{\ast}_4\tau_0,\lll 1,4)
   \stackrel{ R_1^{(3)}}{\longrightarrow}
    (\co_{X_4},\,\co_{X_4}(\ell),\,\lll 1,4,\,\co_{X_4}(-\omega_4-\ell))
       \stackrel{ R_3}{\longrightarrow}\\ \bigskip
   \longrightarrow
    \left(\co_{X_4},\,\co_{X_4}(\ell),\,\co_{X_4}(-\omega_4-\ell),\,
                  \left\{\co_{X_4}(\ell_i-\omega_4-\ell)\right\}_{i=1,2,3,4}\right)
       \stackrel{ L_2}{\longrightarrow}\\
\longrightarrow \renewcommand{\arraystretch}{1.3}
\begin{array}{|c|c|@{}c@{}|}
\hline
   &       &\upl{\co_{X_4}(\ell)}\\
         &       &\co_{X_4}(\ell_1-\omega_4-\ell)\\
\ \quad\co_{X_4}\quad&\quad  F\quad   &\co_{X_4}(\ell_2-\omega_4-\ell)\\
         &       &\co_{X_4}(\ell_3-\omega_4-\ell)\\
         &       &\co_{X_4}(\ell_4-\omega_4-\ell)\\
\hline
\end{array}=\tau\xx4 .
\end{array}
$$
This is the collection corresponding to the solution
$(1,2,1)$ of \x4. To the solution $(1,1,2)$, the collection
$R_1\circ R_2\circ R_2(\tau\xx4)$ corresponds.

The bundle $F$ is obtained as a universal extension
$$
0\longrightarrow\co_{X_4}(-\omega_4-\ell)\longrightarrow F
  \longrightarrow\co_{X_4}(\ell)\longrightarrow 0.
$$
Here, $c_1(F)=-\omega_4,\ r(F)=2$.

In other cases, we present the starting and final collections only and the
sequence of mutations. Verifying details is left to an interested reader.

\smallskip
\x5. Over $X_5$, the mutations $R_1^{(3)}\circ R_1$ of the 4-block collection
$(\LL 4,5(-1),\sigma^{\ast}\tau\xx3)$ result in the desired collection
$$
\settowidth{\plw}{$\ \co_{X_5}(2\ell-\ell_2-\ell_3)\ $}
\tau\xx5=\ \renewcommand{\arraystretch}{1.3}
\begin{array}{|c|c|@{}c@{}|}
\hline
         &       &\upl{\co_{X_5}(-\omega_5)}\\
\shup{\co_{X_5}(\ell_4)}&\shup{\co_{X_5}(\ell)}&\co_{X_5}(2\ell-\ell_1-\ell_2)\\
\shdown{\co_{X_5}(\ell_5)}&\shdown{\co_{X_5}(2\ell-\ell_1-\ell_2-\ell_3)}&
\co_{X_5}(2\ell-\ell_2-\ell_3)\\
         &       &\co_{X_5}(2\ell-\ell_1-\ell_3)\\
\hline
\end{array}\ .
$$

\smallskip
\x{6.1}. The sequence of mutations $R_1^{(3)}\circ L_3\circ R_1$ takes the
4-block collection $(\lll 4,6(-1), \sigma^{\ast}\tau\xx3)$ over $X_6$ to the
desired one,
$$
\settowidth{\plw}{$\ \co_{X_6}(2\ell-\ell_1-\ell_2-\ell_3)\ $}
\tau\xx{6.1}=\renewcommand{\arraystretch}{1.3}
\begin{array}{|c|c|@{}c@{}|}
\hline
\co_{X_6}(\ell_4)&\co_{X_6}(\ell-\ell_1)&
\upl{\co_{X_6}(-\omega_6)\mathstrut}\\
\co_{X_6}(\ell_5)&\co_{X_6}(\ell-\ell_2)&\co_{X_6}(\ell)\\
\co_{X_6}(\ell_6)&\co_{X_6}(\ell-\ell_3)&\co_{X_6}(2\ell-\ell_1-\ell_2-\ell_3)\\
\hline
\end{array}\ .
$$

\smallskip
\x{6.2}. The 4-block collection $(\lll 1,6(-1), \sigma^{\ast}\tau_0)$ over
$X_6$ is taken by
$R_1^{(3)}\circ R_1^{(3)}\circ R_1\circ R_2$ to the desired one,
$$
\settowidth{\plw}{$\ \co_{X_6}(-\omega_6)\ $}
\tau\xx{6.2}=\renewcommand{\arraystretch}{1.3}
\begin{array}{|c|@{}c@{}|c|}
\hline
   &\shdown{\co_{X_6}(\ell)}&
               \co_{X_6}(\ell_1-\omega_6)\quad\co_{X_6}(\ell_2-\omega_6)\\
\quad T_6\quad&\pline{0.3em},{\plw}&\co_{X_6}(\ell_3-\omega_6)\quad\co_{X_6}(\ell_4-\omega_6)\\
   &\shup{\co_{X_6}(-\omega_6)}&
               \co_{X_6}(\ell_5-\omega_6)\quad\co_{X_6}(\ell_6-\omega_6)\\
\hline
\end{array}\ .
$$
Here, $T_6=\sigma^{\ast}_6\,{\rm T}\PP^2(-1))$, $c_1(T_6)=\ell$, $r(T_6)=2$, and
$\mu(T_6)=\ds\frac32$.

\smallskip
\x{7.1}. In this case, the way found by the authors is rather long, and we
divide it into two parts. The 4-block collection
$(\sigma^{\ast}\tau_0, \lll 1,7)$ is taken by
$L_3^{(3)}\circ L_3\circ R_1$ to
$$
\left(\co_{X_7}(\ell+\omega_7),\,\co_{X_7},\,T_7,\,\left\{
\co_{X_7}(\ell-\ell_i)\right\}_{\,i=1,\dots,7}\right) ,
$$
where $T_7=\sigma^{\ast}_7\,{\rm T}\PP^2(-1)$ with the same $r$, $c_1$, and
$\mu$ as those of $T_6$. The latter collection is taken by $R_1^{(3)}\circ
R_1$ to the desired collection
$$
\settowidth{\plw}{$\ \co_{X_7}(\ell-\ell_1),\,\dots\,,\co_{X_7}(\ell-\ell_7)\ $}
\tau\xx{7.1}=
\begin{array}{|c|c|@{}c@{}|}
\hline
         &   &\shdown{\co_{X_7}(-\omega_7)}\\
\quad E_7\quad&\quad T_7\quad&\pline{0.3em},{\plw}\\
         &   &\shup{\co_{X_7}(\ell-\ell_1),\,\dots\,,\co_{X_7}(\ell-\ell_7)}\\
\hline
\end{array}\ .
$$
The bundle $E_7$ is obtained as the extension
$$
0\longrightarrow \co_{X_7}\longrightarrow E_7\longrightarrow
\co_{X_7}(\ell+\omega_7)\longrightarrow 0.
$$
Then $c_1(E_7)=\ell+\omega_7$, $r(E_7)=2$, and $\mu(E_7)=\ds\frac12$.
In the corresponding sequence of mutations, $R_1\circ L_1={\rm id}$, i.e.,
$$
R_1^{(3)}\circ R_1\circ (L_1\circ L_2\circ L_3)\circ L_3=
R_1^{(3)}\circ L_2\circ L_3\circ L_3.
$$

\smallskip
\x{7.2}. The sequence $L_3^{(3)}\circ R_1\circ R_3$ takes
$(\lll 4,7,\sigma^{\ast}\tau\xx3)$ to the 4-block collection
$$\renewcommand{\arraystretch}{1.3}
\left(
  \begin{array}{c} \shup{\co_{X_7}(\ell+\omega_7)}\\
                   \shdown{\co_{X_7}(\omega_7-\ell-\omega_3)}
   \end{array}\ ,\quad\co_{X_7}\ ,\quad
  \begin{array}{c} \co_{X_7}(\ell_4)\\ \co_{X_7}(\ell_5)\\
                    \co_{X_7}(\ell_6)\\ \co_{X_7}(\ell_7)
   \end{array}\ ,\quad
  \begin{array}{c} \co_{X_7}(\ell-\ell_1)\\ \co_{X_7}(\ell-\ell_2)\\
                   \co_{X_7}(\ell-\ell_3)
   \end{array}
\right) ,
$$
which is taken by $R_1^{(3)}\circ R_1$ to the desired collection
$$
\settowidth{\plw}{$\ \co_{X_7}(\ell-\ell_2)\ $}
\tau\xx{7.2}=\renewcommand{\arraystretch}{1.3}
\begin{array}{|c|c|@{}c@{}|}
\hline
&\co_{X_7}(\ell_4)&\upl{\co_{X_7}(-\omega_7)\mathstrut}\\
\quad\shup{E_7}\quad&\co_{X_7}(\ell_5)&\co_{X_7}(\ell-\ell_1)\\
\quad\shdown{E'_7}\quad&\co_{X_7}(\ell_6)&\co_{X_7}(\ell-\ell_2)\\
     &\co_{X_7}(\ell_7)&\co_{X_7}(\ell-\ell_3)\\
\hline
\end{array}\ .
$$
Here $E_7$ is the rank-2 bundle described above and $E'_7$ is the result of
the right shift of the pair $(\co_{X_7}(\omega_7-\ell-\omega_3),
\co_{X_7})$,
$$
0\longrightarrow \co_{X_7}\longrightarrow E'_7\longrightarrow
\co_{X_7}(\omega_7-\ell-\omega_3)\longrightarrow 0.
$$
We have $c_1(E'_7)=\omega_7-\ell-\omega_3=-\ell+\ell_4+\ell_5+\ell_6+\ell_7$, $r(E'_7)=2$, and
$\mu(E'_7)=\ds\frac12$.

\smallskip
\x{7.3}. The 4-block collection $(\co_{\ell_7}(-1),\sigma^{\ast}\tau\xx{6.1})$
over $X_7$ is taken by $R_1^{(3)}\circ R_1$ to the desired one,
$$
\settowidth{\plw}{$\ \co_{X_7}(-\omega_6)=\co_{X_7}(\ell_7-\omega_7)\ $}
\tau\xx{7.3}=\renewcommand{\arraystretch}{1.3}
\begin{array}{|c|c|@{}c@{}|}
\hline
           &                &\shdown{\co_{X_7}(\ell)}\\
           &                &\shdown{\co_{X_7}(2\ell-\ell_1-\ell_2-\ell_3)}\\
   &\shup{\co_{X_7}(\ell-\ell_1)} &\shdown{\co_{X_7}(-\omega_6)
                              =\co_{X_7}(\ell_7-\omega_7)}\\
\quad E''_{7}\quad&\co_{X_7}(\ell-\ell_2)&\pline {0.3em},{\plw}\\
  &\shdown{\co_{X_6}(\ell-\ell_3)}&\shup{\co_{X_7}(\ell_6-\omega_7)}\\
           &                &\shup{\co_{X_7}(\ell_5-\omega_7)}\\
           &                &\shup{\co_{X_7}(\ell_4-\omega_7)}\\
\hline
\end{array}\ .
$$
Here, the bundle $E''_{7}$ is the result of the right shift of the torsion
sheaf $\co_{\ell_7}(-1)$ over the block
$(\co_{X_7}(\ell_4),\co_{X_7}(\ell_5),\co_{X_7}(\ell_6))$,
$$
0\longrightarrow (\co_{X_7}(\ell_4)\oplus\co_{X_7}(\ell_5)\oplus\co_{X_7}(\ell_6))
\longrightarrow E''_{7}\longrightarrow
\co_{\ell_7}(-1)\longrightarrow 0.
$$
We have $c_1(E''_{7})=\ell_4+\ell_5+\ell_6+\ell_7$, $r(E''_{7})=3$, and
$\mu(E''_{7})=\ds\frac{4}{3}$.

\smallskip
\x{8.1}. Over $X_8$, the collection $(\sigma^{\ast}(\tau_0(-1)),\lll 1,8)$ is
taken by $L_3^{(3)}\circ L_3\circ R_1$ to the 4-block collection
$$
\left(\co_{X_8}(-\omega_8) ,\, \co_{X_8}(-\ell) ,\, \sigma^{\ast}(\TPP(-2))
,\, \left\{\co_{X_8}(-\ell_i)\right\}_{i=1,\dots,8}\right).
$$
Applying $L_1\circ L_1$ to the latter collection, we obtain the desired one,
$$
\settowidth{\plw}{$\ \co_{X_8}(-\ell_1),\,\dots\,,\co_{X_8}(-\ell_8)\ $}
\tau\xx{8.1}=
\begin{array}{|c|c|@{}c@{}|}
  \hline
       &          &\shdown{\co_{X_8}(-\omega_8)}\\
\quad E_8\quad&\quad F_8\quad&\pline {0.3em},{\plw}\\
       &          &\shup{\co_{X_8}(-\ell_1),\,\dots\,,\co_{X_8}(-\ell_8)}\\
  \hline
\end{array}\ .
$$
Here $E_8=L_{\co_{X_8}(-\omega_8)}\co_{X_8}(-\ell)$ and
$F_8=L_{\co_{X_8}(-\omega_8)}\sigma^{\ast}(\TPP(-2))$. Both mutations are
extensions,
$c_1(E_8)=-2\omega_8-\ell$, $c_1(F_8)=-\omega_8-\ell$,
$r(E_8)=r(F_8)=3$, $\mu(E_8)=-\ds\frac{5}{3}$, and
$\mu(F_8)=-\ds\frac{4}{3}$.

\smallskip
\x{8.2}. Applying $L_3\circ L_3$ to $(\lll 4,8(-1),\sigma^{\ast}\tau\xx3)$,
we obtain the 4-block collection
$$\renewcommand{\arraystretch}{1.3}
\left(
\lll 4,8(-1)\, ,\ \co_{X_8}\, ,\
\begin{array}{c} T_8\\ T'_8 \end{array} ,\
\begin{array}{c} \co_{X_8}(\ell-\ell_1)\\
\co_{X_8}(\ell-\ell_2)\\ \co_{X_8}(\ell-\ell_3)
\end{array}
\right) .
$$
Here, the bundle $T_8=\sigma^{\ast}\TPP(-1)$ is the left shift of
$\co_{X_8}(2\ell-\ell_1-\ell_2-\ell_3)$ over the block
$\left\{\co_{X_8}(\ell-\ell_i)\right\}_{i=1,2,3}$. This is the unique (according to
1.3) exceptional bundle of rank 2 over $X_8$ with $c_1=\ell$. The bundle
$T'_8$ is the left shift of $\co_{X_8}(\ell)$ over the same block. Then
$c_1(T'_8)=-\omega_3-\ell$, $r(T'_8)=2$, and $\mu(T'_8)=\mu(T_8)=\ds\frac{3}{2}$.

Applying $R_1^{(3)}\circ R_1\circ R_1$ to the latter 4-block collection
results in the desired collection
$$
\settowidth{\plw}{$\ \co_{X_8}(\ell_4-\omega_8)\ %
                  \co_{X_8}(\ell_5-\omega_8)\ \co_{X_8}(\ell_6-\omega_8)\ $}
\tau\xx{8.2}=\renewcommand{\arraystretch}{1.3}
\begin{array}{|c|c|@{}c@{}|}
\hline
    &\quad\shdown{T_8}\quad &\ \co_{X_8}(\ell_4-\omega_8)\ %
                  \co_{X_8}(\ell_5-\omega_8)\ \co_{X_8}(\ell_6-\omega_8)\ \\
\quad E'_8\quad&  &\upl{\co_{X_8}(\ell_7-\omega_8)\ \co_{X_8}(\ell_8-\omega_8)}\\
    &\shup{T'_8} & \co_{X_8}(\ell-\ell_1)\ \co_{X_8}(\ell-\ell_2)\ \co_{X_8}(\ell-\ell_3)\\
\hline
\end{array}\ .
$$
The bundle $E'_8$ is the right shift of $\co_{X_8}$ over
$\left\{\co_{X_8}(\ell_i)\right\}_{i=4,\dots,8}$. Computation shows that
$c_1(E'_8)=\sum\limits_{i=4}^8 \ell_i$, $r(E'_8)=4$, and
$\mu(E'_8)=1\frac{1}{4}$.

\smallskip
\x{8.3}. Over $X_8$, the 4-block collection
$(\LL 7,8(-1),\sigma^{\ast}\tau\xx{6.1})$ is taken by
$R_1^{(3)}\circ L_3\circ R_1$ to the desired collection
$$
\settowidth{\plw}{$\ \co_{X_8}(\ell_4-\omega_8)\ %
                  \co_{X_8}(\ell_5-\omega_8)\ \co_{X_8}(\ell_6-\omega_8)\ $}
\tau\xx{8.3}=\renewcommand{\arraystretch}{1.3}
\begin{array}{|c|c|@{}c@{}|}
\hline
\shdown{\sigma^{\ast}E''_7}&T_8&\shdown{\co_{X_8}(\ell_4-\omega_8)\ %
                  \co_{X_8}(\ell_5-\omega_8)\ \co_{X_8}(\ell_6-\omega_8)\ }\\
         &\quad T'_8\quad &\pline {0.3em},{\plw}\\
\quad\shup{E''_8}\quad&T''_8&\shup{\co_{X_8}(\ell-\ell_1)\ \co_{X_8}(\ell-\ell_2)\ %
                          \co_{X_8}(\ell-\ell_3)}\\
\hline
\end{array}\ .
$$
The bundle $\sigma^{\ast}E''_7$ is the right shift of the torsion sheaf
$\co_{\ell_7}(-1)$ over $\left\{\co_{X_8}(\ell_i)\right\}_{i=4,5,6}$. This mutation
is the lifting under $\sigma\colon X_8\rightarrow X_7$ of the mutation
described in item \x{7.3}. The bundle $E''_8$ is the right shift of
$\co_{\ell_8}(-1)$ over the same block
$\left\{\co_{X_8}(\ell_i)\right\}_{i=4,5,6}$. The correspondent mutation is an
extension, and we have
$c_1(E''_8)=\ell_4+\ell_5+\ell_6+\ell_8$, $r(E'_8)=4$, and
$\mu(E''_8)=\mu(\sigma^{\ast}E''_7)=\ds\frac{4}{3}$.

The bundles $T_8$ and $T'_8$ are described in the previous item,
and $T''_8$ is the left shift of $\co_{X_8}(-\omega_6)$ over
$\left\{\co_{X_8}(\ell-\ell_i)\right\}_{i=1,2,3}$ (the middle block in
$\tau\xx{6.1}$). We have $c_1(T''_8)=\omega_6-\omega_3$, $r(T''_8)=2$, and
$\mu(T''_8)=\mu(T'_8)=\mu(T_8)=\ds\frac{3}{2}$.

\smallskip
\x{8.4}. Applying $R_1^{(3)}\circ L_3\circ R_1\circ L_2$, to
$(\sigma^{\ast}\tau\xx3,\lll 4,8)$, we obtain the collection corresponding to
the solution $(5,2,1)$,
$$
\settowidth{\plw}{$\ \co_{X_8}(2\ell-\ell_1-\ell_2-\ell_3)\ $}
\tau\xx{8.4}=\renewcommand{\arraystretch}{1.3}
\begin{array}{|c|c|@{}c@{}|}
\hline
           &\quad F_{48}\quad   &\co_{X_8}(\ell)\\
           &      F_{58}      &\upl{\co_{X_8}(2\ell-\ell_1-\ell_2-\ell_3)}\\
\quad E'''_8\quad & F_{68}   &\co_{X_8}(\ell-\ell_1-\omega_8)\\
           & F_{78}      &\co_{X_8}(\ell-\ell_2-\omega_8)\\
           & F_{88} &\co_{X_8}(\ell-\ell_3-\omega_8)\\
\hline
\end{array}\ .
$$
To the solution $(5,1,2)$, the collection
$L_2\circ L_1\circ L_1(\tau\xx{8.4})$ corresponds.

Here, $E'''_8$ is the right shift of $\co_{X_8}$ over
$\left\{\co_{X_8}(\ell-\ell_i)\right\}_{i=1,2,3}$. We have
$c_1(E'''_8)=-2\omega_3$, $r(E'''_8)=5$, and $\mu(E'''_8)=2\frac{2}{5}$.
The bundles $F_{i8}$ are obtained as the shifts of the torsion sheaves
$\co_{\ell_i}$ over
$\left(\co_{X_8}(\ell),\co_{X_8}(2\ell-\ell_1-\ell_2-\ell_3)\right)$. We have
$c_1(F_{i8})=-\omega_3-\ell_i$, $i=4,5,6,7,8$, $r(F_{i8})=2$, and
$\mu(F_{i8})=2\frac{1}{2}$.
This completes the proof of the proposition.

\smallskip
Note that the lifting procedure used to obtain all minimum 3-block
collections in the latter proof, is revertible. The reverse procedure is the
following: a distinguished block is divided into two blocks, and to the
obtained 4-block collection, the inverse sequence of mutations is applied.
Moreover, the same procedure can be applied to an arbitrary minimum 3-block
collection $({\cal E},{\cal F},{\cal G})$ over a Del Pezzo surface as well
since the invariants $(r({\cal E}),r({\cal F}),r({\cal G}))$,
$(\alpha,\beta,\gamma)$, and
$(\chi({\cal E},{\cal F}),\chi({\cal F},{\cal G}))$ are determined by the
equation itself and, therefore, take the same values as for one of the
minimum collections obtained in the latter proof (see also (16), (17)).
Here, it is clear which block is distinguished, and into what parts it should
be divided. One easily sees that the ranks of sheaves contained in an
exceptional collection and the dimensions of $\ext$ spaces defined by
(6) can be uniquely determined for the collection obtained by a mutation.
Thus, applying the above-mentioned procedure to
$({\cal E},{\cal F},{\cal G})$, we get a 4-block collection, one of whose
blocks consists of zero-rank sheaves. Taking into account 1.6, we arrive at
the statement below.

\smallskip
{\bf 4.3. Proposition. }{\it Let $S$ be a Del Pezzo surface other than\/
$\PP^2$ or\/ $\PP^1\times\PP^1,$ $\tau$ be a minimum\/ $3$-block collection
of type\/ $(\alpha,\beta,\gamma)$, $\alpha\le\beta\le\gamma$, over\/
$S$. Then a divisor $D$ over $S$ and a sequence of mutations\/ $($which does
not preserve the\/ $3$-block structure$)$ exist such that this sequence takes\/
$\tau(D)$ to\/ $(\sigma^{\ast}\tau',(\co_{e_1},\dots,\co_{e_m}))$, where\/
$e_1 ,\dots,e_m$ are pairwise nonintersecting exceptional curves,\/
$\sigma\colon S\rightarrow S'$ is the monoidal transform with center\/
$\{e_1 ,\dots,e_m\}$, and\/ $\tau'$ is a minimum\/ $3$-block collection
over\/ $S'$.}

\bigskip
\begin{center}{\large\bf 5. The action of the Weyl group on complete
3-block collections}\end{center}

\medskip
{\bf 5.1. }In this section, we denote by $X_r$ a Del Pezzo surface of degree
$9-r$, i.e., we do not fix an identification of $X_r$ with a plane with
$r$ blown-up points. Denote by $I_r\subset{\rm Pic}\,X_r$ the set of classes
of exceptional curves, the latter being characterized by the equalities
$e^2=e\cdot\omega_r=-1$. Let ${\cal R}_r\subset{\rm Pic}\,X_r$ be the set of
vectors $s$ such that $s^2=-2$ and $s\cdot\omega_r=0$. The form obtained by
changing sign of the intersection form on
${\rm Pic}\,X_r\otimes_{\Bbb Z}{\Bbb R}$ induces the structure of a Euclidean
space on the orthogonal complement to $\omega_r$. The set ${\cal R}_r$ is a
root system in it. The Weyl group  $W({\cal R}_r)$ of this system, which is
generated by symmetries with respect to the roots, coincides with the group
of automorphisms of the lattice ${\rm Pic}\,X_r$ which preserve $\omega_r$ and
the intersection form and also with the group of permutations of elements
of $I_r$ which preserve pairwise intersection indices of the elements (see
[14]). We denote this group by $W_r$ and call it the Weyl group. Over
$\PP^1\times\PP^1$, there is only one (up to the sign) divisor with square
$-2$, i.e., the Weyl group for the quadric is $\ZZ_2$.

\smallskip
{\bf Agreement.}
{\it All statements of this section concerning $X_r$ and $W_r$, except
Proposition\/ $5.4$, concern\/ $\PP^1\times\PP^1$ as well.}

\smallskip
Fix the isomorphism
$$
v:\:{\rm K}_0(X_r)\longrightarrow \ZZ\oplus\pic X_r\oplus\ZZ,\quad
E\stackrel{v}{\longmapsto} (r(E),c_1(E),2\ch_2(E)),
$$
where $\ch_2(E)=c_1^2(E)/2-c_2(E)$ is the second component of the Chern
character. The action of the Weyl group on $\pic X_r$ is in a natural way
extended to the action on K$_0(X_r)$, $r$ and $2\ch_2 $ do not change
therewith. We want to define the action of $W_r$ on the set of exceptional
collections in concordance with the action on K$_0(X_r)$. Put
$$
g\co_{X_r}(D)\triangleq\co_{X_r}(gD)\qquad {\rm and}\qquad
g\co_{\ell}(m)\triangleq\co_{g\ell}(m).
$$

\smallskip
{\bf 5.2. Lemma. }{\it Let\/ $(E,F)$ and\/ $(E',F')$ be exceptional pairs of
sheaves, where\/ $v(E')=gv(E)$ and $v(F')=gv(F)$. Then the statements below
hold,}

(a) $v(R_F'E')=gv(R_FE)$ {\it and\/} $v(L_E'F')=gv(L_EF)$;

(b) $\dim\ext^i(E',F')= \dim\ext^i(E,F),\ \forall i$.

\smallskip
{\sc Proof}. A pair $(E,F)$ is the elementary two-block collection.
Propositions 2.2 and 2.3 imply that the mutation types of this pair depend on
$r(E)$, $r(F)$, $c_1(E)\cdot K$, and $c_1(F)\cdot K$ only (the value
of $\chi (E,F)=\chiminus (E,F)$ is determined by (6)). The action of the Weyl
group preserves the intersection form and the canonical class. Hence, the
above-mentioned invariants for $(E',F')$ take the same values as for
$(E,F)$. Thus, both pairs have the same type of left and right mutations.
Now, (a) is verified by direct computation using the exact sequences of
Sec.~1.7, additivity of $v$, and linearity of $g$. The validity of (b)
follows from the equality $\chi (E,F)=\chi (E',F')$ and the classification
of exceptional pairs (see Sec.~1.7). The lemma is proved.

\smallskip
{\bf 5.3. Proposition--definition. }{\it Let\/ $\tau = (E_1,\dots ,E_n)$ be
an exceptional collection and\/ $g\in W_r$. Then a unique exceptional
collection\/ $\tau' = (E'_1,\dots ,E'_n)$ exists such that\/
$v(E'_i)=gv(E_i)$ for\/ $i=1,\dots ,n$.

Put\/ $g\tau\triangleq\tau'$ and\/ $gE_i\triangleq E'_i$}.

\smallskip
{\sc Proof}. According to [11, 6.11], $\tau$ is included in a complete
exceptional collection. Therefore, without loss of generality, we assume
$\tau$ to be complete. By the constructivity theorem [11,~7.7], $\tau$ is
obtained by a sequence of mutations (in the sense of Sec.\ 1.10) from the
complete collection
$$
\tau_1=\left(\co_{\ell_1}(-1),\dots ,\co_{\ell_r}(-1),
\co_{X_r},\co_{X_r}(\ell_0),\co_{X_r}(2\ell_0)\right),
$$
for which the proposition obviously holds. Applying this sequence
of mutations to the
complete exceptional collection $g\tau_1$, we obtain the desired exceptional
collection according to the preceding lemma. The uniqueness follows from
Proposition 1.3.

\smallskip
Thus, the Weyl group acts on the set of exceptional collections and preserves
the ranks of sheaves and the dimensions of $\ext$ spaces. According to 5.2a,
this action commutes with mutations. In what follows, we are interested in
the action of $W_r$ on complete 3-block collections of sheaves.

\smallskip
{\bf 5.4. Proposition. }{\it Let\/ $(\EE,\FF,\GG)$ and\/
$(\EE',\FF',\GG')$ be two minimum collections of type\/
$(\alpha,\beta,\gamma)$, $\alpha\le\beta\le\gamma$, corresponding to the same
minimum solution\/ $($the latter is essential for the equations of group\/
{\rm IV)}. Then\/ an element $g\in W_r$ and a divisor\/ $D\in{\rm Pic}\,X_r$
exist such that}
$$
(\EE,\FF,\GG)=g\big(\EE'(D),\FF'(D),\GG'(D)\big).
$$

\smallskip
{\sc The proof} is carried out by induction on $r$. For $r=0$, the statement
is trivial due to triviality of $W_0$. Let the proposition hold for all
$p<r$. By 4.3, there exist a divisor $D_1\in{\rm Pic}\,X_r$ and a sequence of
mutations $\Phi$ (which does not preserve the 3-block structure) that takes
$\big(\EE(D_1),\FF(D_1),\GG(D_1)\big)$ to a 4-block collection of the form
$(\sigma^*\tau,{\cal L}_e)$. Here $e=\{e_1,\dots,e_{r-p}\}$ is a set of
$(-1)$-curves with $e_i\cdot e_j=0$ for $i\ne j$, the block ${\cal L}_e$
consists of the sheaves $\OO_{e_i}$, $i=1,\dots,r-p$, the morphism
$\sigma:X_r\to X_p$ is the blowing down for $e$, and $\tau$ is a 3-block
collection over $X_p$. The same sequence of mutations $\Phi$ takes
$\big(\EE'(D_2),\FF'(D_2),\GG'(D_2)\big)$ for some $D_2\in{\rm Pic}\,X_r$ to
an analogous collection $(\sigma'^*\tau',{\cal L}_{e'})$, where
$e'=\{e'_1,\dots,e'_{r-p}\}$ is another set of pairwise nonintersecting
$(-1)$-curves, $\sigma':X_r\to X_p$ is the blowing down for $e'$, and $\tau'$
is a minimum collection over $X_p$ of the same type as $\tau$.

Consider an element $g_0\in W_r$ that takes $e'_i$ to $e_i$,
$i=1,\dots,r-p$. Then
$$
g_0(\sigma'^*\tau',{\cal L}_{e'})=(\sigma^*\tau'',{\cal L}_e).
\eqno(19)
$$
Here, $\tau''$ is the minimum collection over $X_p$ of the same type as
$\tau$. By the induction assumption, an element $h\in W_p$ and a divisor
$D_0\in{\rm Pic}\,X_p$ exist such that $\tau=h(\tau''(D_0))$. Identifying
$W_p$ with a subgroup in $W_r$ that preserves all $(-1)$-curves $e_i$,
$i=1,\dots,r-p$, and identifying ${\rm Pic}\,X_p$ with the orthogonal
complement to $e$ in ${\rm Pic}\,X_r$ (with respect to the intersection
form), we obtain
$$
(\sigma^*\tau,{\cal L}_e)=h(\sigma^*\tau''(D_0),{\cal L}_e)=
h(\sigma^*\tau''(D_0),{\cal L}_e(D_0))
\eqno(20)
$$
(note that $\OO_{e_i}(D_0)=\OO_{e_i}$ since $e_i\cdot D_0=0$).

Combining (19) and (20) together, we get
$$
(\sigma^*\tau,{\cal L}_e)=h\big((g_0\sigma'^*\tau')(D_0),{\cal L}_e(D_0)\big)
=hg_0\big(\sigma'^*\tau'(g_0^{-1}D_0),{\cal L}_{e'}(g_0^{-1}D_0)\big).
$$

Applying to both sides of this equality the sequence of mutations
$\Phi^{-1}$, and taking into account that mutations commute with twisting
and the action of $W_r$, we obtain
$$
\big(\EE(D_1),\FF(D_1),\GG(D_1)\big)=hg_0
\big(\EE'(D_2+g_0^{-1}D_0),\FF'(D_2+g_0^{-1}D_0),\GG'(D_2+g_0^{-1}D_0)\big).
$$
Hence, the desired statement easily follows; it suffices to put
$g=hg_0$ and $D=D_2+g_0^{-1}D_0+g^{-1}D_1$.

\medskip
{\bf 5.5. }On the set of complete 3-block collection,

(1) the group ${\rm Pic}\,X_r$ acts by tensoring by invertible sheaves;

(2) the 3-string braid group $B(3)$ acts by mutations;

(3) the Weyl group $W_r$ acts.\\
The action of the braid group commutes with the action of ${\rm Pic}\,X_r$ and
$W_r$.

We call two 3-block collections {\it equivalent}, if one is obtained from
another by tensoring. The preceding proposition means that $W_r$ acts
transitively on the equivalence classes of minimum collections of type
$(\alpha,\beta,\gamma)$, where $\alpha\le\beta\le\gamma$. By Theorem 3.10,
any orbit with respect to $B(3)$ contains minimum collections of this
type. Hence, $W_r$ acts transitively on $(B(3)\times{\rm Pic}\,X_r)$-orbits
of complete 3-block collections of a given structure.

For brevity, we call the set of all collections obtained from a complete
3-block collection by mutations and tensoring (i.e., its orbit under the
action of $B(3)\times{\rm Pic}\,X_r$) the {\it orbit} of this
collection. The number of orbits of complete 3-block collections of a given
structure is finite since $W_r$ is finite. We are going to compute these
numbers. Note that a structure $\{\alpha,\beta,\gamma\}$ of a complete
3-block collection corresponds to a Markov-type equation in the table of
Sec.\ 3.5.

\medskip
{\bf 5.6. Proposition. }{\it The number\/ $C$ of equivalence classes of
minimum collections of type\/ $(\alpha,\beta,\gamma)$,
$\alpha\le\beta\le\gamma$, which lie in one orbit is finite. These numbers
are given in the bottom row of the table below (in the top row, equation
numbers are given).}
{\addtolength{\tabcolsep}{-1.1pt}\begin{center}
\renewcommand{\arraystretch}{1.3}
\begin{tabular}{|c|c|c|c|c|c|c|c|c|c|c|c|c|c|}
\hline
\x1&\x2&\x3&\x4&\x5&\x{6.1}&\x{6.2}&\x{7.1}&\x{7.2}&\x{7.3}%
&\x{8.1}&\x{8.2}&\x{8.3}&\x{8.4}\\
\hline
1 & 1 & 1 & 2 & 2 & 3 & 1 & 2 & 2 & 1 & 1 & 1 & 1 & 2 \\
\hline
\end{tabular}
\end{center}}

\smallskip
{\sc Proof}. For the equation \x1, this follows from the well-known
fact that any exceptional collection over $\PP^2$ which consists of invertible
sheaves has the form
$({\cal O}_{\PP^2}(m-1),{\cal O}_{\PP^2}(m),{\cal O}_{\PP^2}(m+1))$.
For \x2, the statement is also well-known and follows from the description of
exceptional collections over $\PP^1\times\PP^1$ that consist of invertible
sheaves (see [5,~5.6]). Let us consider the remaining equations.

The main tool in the proof is mutations of 3-block helices; see 1.11.
A helix $[\EE,\FF,\GG]$ generated by a complete 3-block collection
$(\EE,\FF,\GG)$ is an infinite sequence of blocks
$$
\big(\dots,\:\EE(K),\:\FF(K),\:\GG(K),\:\EE,\:\FF,\:\GG,\:
\EE(-K),\:\FF(-K),\:\GG(-K),\:\dots\big).
$$
Any three successive blocks of this sequence form a complete 3-block
exceptional collection called a foundation of the helix. Evidently, a helix
is uniquely determined by any foundation. Mutations of foundations in the
sense of Sec.\ 2 define {\it mutations of helices}. For example,
the mappings $[\EE,\FF,\GG]\mapsto [L_{\EE}\FF,\EE,\GG]$ and
$[\EE,\FF,\GG]\mapsto [\EE,\GG,R_{\GG}\FF]$ are helix mutations which
coincide since $L_{\EE}\FF=L_{\EE}L_{\GG}R_{\GG}\FF=R_{\GG}\FF$. The ranks
and lengths in any foundation of a helix satisfy one of the equations (3),
and mutations of a helix induce mutations of the solutions of this equation,
i.e., a route along the corresponding graph of solutions.

\medskip
{\bf Lemma. }{\it Assume that a sequence of mutations\/ $\Phi$ of a\/
$3$-block helix\/ $[\EE,\FF,\GG]$ induces a cyclic route along the solution
graph, the route containing no loop at the minimum solution. Then\/
$\Phi$ preserves the\/ $3$-block helix.}

{\sc Proof}. Let $\Psi$ be the sequence of solution mutations induced by
$\Phi$, and a solution $w$ be the starting and final point of the route.
The solution graph contains no cycles; hence, $\Psi=M\Psi_1M$. Here, $M$ is
one of the mutations $M_x,M_y,M_z$, and $\Psi_1$ is a sequence of solution
mutations inducing a cyclic route that starts and finishes at $M(w)$. Using
induction on the number of mutations in $\Psi$, we may consider that the
sequence of helix mutations which induces $\Psi_1$ preserves the helix
$\widetilde{M}[\EE,\FF,\GG]$, where $\widetilde{M}$ induces $M$.
Hence, $\Phi[\EE,\FF,\GG]=\widetilde{M}^2[\EE,\FF,\GG]$.
The square of a solution mutation over one of the variables can only be
induced by one of the mutations ${\rm id}=R_1L_1=L_1R_1=R_2L_2=L_2R_2$, or
$R_2R_1$, or $L_1L_2$ of a suitable foundation of the helix. We have
$R_2R_1(\EE',\FF',\GG')=(\FF',\GG',\EE'(-K))$ and
$L_1L_2(\EE',\FF',\GG')=(\GG'(K),\EE',\FF')$; thus, all mentioned mutations
of a foundation preserve the helix. The lemma is proved.

\smallskip
{\it Remark}. The lemma is also valid under the assumption that a cyclic
route along the solution graph contains an even number of loops at the
minimal solution. This can easily be checked using the lemma.

\smallskip
We return to the proof of the proposition. Let $(\EE,\FF,\GG)$ and
$(\EE',\FF',\GG')$ be minimum collections of type $(\alpha,\beta,\gamma)$,
$\alpha\le\beta\le\gamma$ (corresponding to the same minimum solution, for
the case of equations of group IV), and $\Phi$ be a sequence of mutations
that takes the first collection to the second one. Consider the sequence
$\Psi$ of mutations of solutions of the corresponding Markov-type equation
which is induced by $\Phi$. It defines a cyclic route along the solution
graph with starting and final point at the minimum solution. Two cases are
possible.

(a) If $\Psi$ contains an even number of loops, then, as is proved above,
the collections $(\EE,\FF,\GG)$ and $(\EE',\FF',\GG')$ are foundations of the
same helix. For all equations except \x{6.1}, $\alpha <\beta$ or
$\beta <\gamma$, and coincidence of types of these collections implies
their equivalence, i.e., $(\EE,\FF,\GG)=(\EE'(mK),\FF'(mK),\GG'(mK))$.

(b) If $\Psi$ contains an odd number of loops, then
$\Psi=\Psi_1M\Psi_2$, where $M$ is the loop and $\Psi_1,\ \Psi_2$ are the
routes that start and finish at the minimum solution and contain even
numbers of loops. Then $\Phi=\Phi_1\widetilde{M}\Phi_2$, where $\Phi_1$ and
$\Phi_2$ are sequences of mutations of 3-block collections which induce
$\Psi_1$ and $\Psi_2$, and $\widetilde{M}$ induces $M$. As is proved above,
$\Phi_1$ and $\Phi_2$ preserve the 3-block helix. Hence,
$$
[\EE',\FF',\GG']=\Phi[\EE,\FF,\GG]=\widetilde{M}[\EE,\FF,\GG].
$$

Let us consider each group of equations separately.

\smallskip
I. The solution graph is loopless; hence, the case (a) is possible only.

For the equation \x{6.1}, $\alpha=\beta=\gamma=3$, and one can consider (up
to tensoring by $mK_{X_6}$) that the collection $(\EE',\FF',\GG')$ 
coincides with one of the
collections $(\EE,\FF,\GG)$, $(\GG(K),\EE,\FF)$, or $(\FF(K),\GG(K),\EE)$. It
is easy to check that these collections are not equivalent for
$(\EE,\FF,\GG)=\tau_{\xx{6.1}}$. Hence, by 5.4, they are not equivalent in a
general case as well.

For \x{8.1}, the statement is proved in (a).

\smallskip
II. Consider the equation \x5. In the case (a), collections
$(\EE,\FF,\GG)$ and $(\EE',\FF',\GG')$ are equivalent. In the case (b),
applying the action of the Weyl group and tensoring the second collection by
$mK_{X_5}$, we can obtain $(\EE,\FF,\GG)=\tau_{\xx5}$ and
$(\EE',\FF',\GG')=R_1R_2R_2\tau_{\xx5}$. Direct computations show that the
collections $\tau_{\xx5}$ and $R_1R_2R_2\tau_{\xx5}$ are not equivalent.
Hence, in the case (b), the collections $(\EE,\FF,\GG)$ and
$(\EE',\FF',\GG')$ are not equivalent too.

For \x{7.1} and \x{7.2}, the reasoning is similar.

For \x{8.2}, we have $\alpha<\beta<\gamma$. The types of foundations obtained
by a single mutation from $[\EE,\FF,\GG]$ are odd permutations of
$(\alpha,\beta,\gamma)$. Hence, the case (b) is impossible.

\smallskip
III. For all equations of this group, we have $\alpha<\beta<\gamma$, and the
reasoning is similar to that for \x{8.2} given above.

\smallskip
IV. For equations of this group, different minimum solutions correspond to
the helices $[\EE,\FF,\GG]$ and $\widetilde{M}[\EE,\FF,\GG]$, so the case (b)
is also impossible here. But the number of equivalence classes of the minimum
collections under consideration equals two, as well as the number of the
minimum solutions. This completes the proof of the proposition.

\medskip
{\bf 5.7. The number of orbits. } In this final subsection, we describe the
computation of the number of orbits under the action of
$B(3)\times{\rm Pic}\,X_r$ on complete 3-block collections of a given
structure. For $\PP^2$ and $\PP^1\times \PP^1$, there is one orbit; this is
shown in [4] and [19]. In the other cases, we find at first the number $N$ of
equivalence classes of minimum collections of type $(\alpha,\beta,\gamma)$,
where $\alpha\le\beta\le\gamma$. By 4.3, such a collection $\tau$ over a
Del Pezzo surface $S$ with the help of tensoring and a sequence of mutations
(determined by the type of the collection from the proof of Proposition 4.2)
can be reduced to the form
$(\sigma^{\ast}\tau',(\co_{e_1},\dots,\co_{e_m}))$, where
$e_1,\dots,e_m$ are pairwise nonintersecting exceptional curves,
$\sigma\colon S\longrightarrow S'$ is the monoidal transform with center
$\{e_1 ,\dots,e_m\}$, and $\tau'$ is a minimum 3-block collection over
$S'$. Note that this sequence of mutations does not preserve the 3-block
structure; and before applying it, one should find the distinguished block
and divide it into two blocks, the length of one of them being $m$.
Thus, we have
$$
N\cdot {n\choose m}=N'\cdot\Big(\mbox{the number of sets }\,
 \{e_1 ,\dots,e_m\},\ e_i\cdot e_j=0 \Big),
$$
where $n$ is the length of the distinguished block and $N'$ is the number of
equivalence classes of minimum collections over $S'$ with the same type as
that of $\tau'$. Computing the values of $N$ for all Markov-type equations
with the help of this formula and dividing them by the correspondent values
of $C$ from Proposition 5.6, we obtain the number of orbits. The results are
presented in the table below.
\vskip \abovedisplayskip
\begin{figure}
\begin{center}\renewcommand{\arraystretch}{1.3}
\begin{tabular}{|c|c|c|r|c|r|} \hline
Equation  & Surface     & $(\alpha,\beta,\gamma)$  &
                                       $ N\quad$&$\quad C\quad$& Number\\
number    &             &            &             &   & of orbits\\
\hline
\x1       & $\PP^2$     &  $(1,1,1)$ &     1$\quad$& 1 &     1    \\
\x2       & $\PP^1
           \times\PP^1$ &  $(1,1,2)$ &     1$\quad$& 1 &     1    \\
\x3       & $X_3$       &  $(1,2,3)$ &     1$\quad$& 1 &     1    \\
\x4       & $X_4$       &  $(1,1,5)$ &     2$\quad$& 2 &     1    \\
\x5       & $X_5$       &  $(2,2,4)$ &    20$\quad$& 2 &    10    \\
\x{6.1}   & $X_6$       &  $(3,3,3)$ &   240$\quad$& 3 &    80    \\
\x{6.2}   & $X_6$       &  $(1,2,6)$ &    36$\quad$& 1 &    36    \\
\x{7.1}   & $X_7$       &  $(1,1,8)$ &    72$\quad$& 2 &    36    \\
\x{7.2}   & $X_7$       &  $(2,4,4)$ &  2520$\quad$& 2 &  1260    \\
\x{7.3}   & $X_7$       &  $(1,3,6)$ &   672$\quad$& 1 &   672    \\
\x{8.1}   & $X_8$       &  $(1,1,9)$ &  1920$\quad$& 1 &  1920    \\
\x{8.2}   & $X_8$       &  $(1,2,8)$ &  8640$\quad$& 1 &  8640    \\
\x{8.3}   & $X_8$       &  $(2,3,6)$ & 80640$\quad$& 1 & 80640    \\
\x{8.4}   & $X_8$       &  $(1,5,5)$ & 96768$\quad$& 2 & 48384    \\
\hline
\end{tabular}
\end{center}
\end{figure}
\vskip \abovedisplayskip

\newpage
\begin{center}{\large\bf REFERENCES}\end{center}

\medskip
\par\noindent\hangindent=\parindent
$\phantom{1}$%
1.~M.~F.~Atiyah, ``Vector bundles over an elliptic curve,'' in: {\it Proc.\
Lond.\ Math.\ Soc.}, {\bf7}, (1957), pp.~414--452.
\par\noindent\hangindent=\parindent
$\phantom{1}$%
2.~A.~I.~Bondal, ``Representations of associative algebras and coherent
sheaves,'' {\it Math. USSR Izvestiya}, {\bf34} (1990), No.~1, 23--42.
\par\noindent\hangindent=\parindent
$\phantom{1}$%
3.~A.~L.~Gorodentsev, ``Transformations of exceptional bundles on
 ${\Bbb P}^n$,'' {\it Math. USSR Izvestiya}, {\bf32} (1989), No.~1, 1--14.
\par\noindent\hangindent=\parindent
$\phantom{1}$%
4.~A.~L.~Gorodentsev and A. N. Rudakov, ``Exceptional vector bundles on
the projective spaces,'' {\it Duke Math.\ J.}, {\bf54}, No.~1, 115--130
(1987).
\par\noindent\hangindent=\parindent
$\phantom{1}$%
5.~A.~L.~Gorodentsev, ``Exceptional bundles over surfaces with a moving
anticanonical class,'' {\it Math. USSR Izvestiya}, {\bf33} (1989), No.~1,
67--83.
\par\noindent\hangindent=\parindent
$\phantom{1}$%
6.~J.-M.~Drezet and J. Le Potier, ``Fibr\'es stables et fibr\'es
exceptionnels sur ${\Bbb P}_2$,'' {\it Ann.\ Scient.\ ENS}, {\bf18}, No.~2,
193--244 (1985).
\par\noindent\hangindent=\parindent
$\phantom{1}$%
7.~J.-M.~Drezet, ``Fibr\'es exceptionnels et suite spectrale de Beilinson
g\'en\'eralis\'ee sur ${\Bbb P}_2({\Bbb C})$,'' {\it Ann.\ Math.},
{\bf 275}, 25--48 (1986).
\par\noindent\hangindent=\parindent
$\phantom{1}$%
8.~J.-M.~Drezet,. ``Fibr\'es exceptionnels et vari\'et\'es de modules de
faisceaux semi-stables sur ${\Bbb P}_2({\Bbb C})$,'' {\it J.\ Reine Angew.\
Math.}, {\bf 380}, 14--58 (1987).
\par\noindent\hangindent=\parindent
$\phantom{1}$%
9.~J.-M.~Drezet, ``Vari\'et\'es de modules extr\'emales de faisceaux
semi-stables sur ${\Bbb P}_2({\Bbb C})$,'' {\it Ann.\ Math.}, {\bf290},
No.~4, 727--770 (1991).
\par\noindent\hangindent=\parindent
10.~S.~Yu.~Zyuzina, ``Constructibility of exceptional pairs of vector bundles
on a quadric,'' {\it Math. USSR Izvestiya}, {\bf42} (1993), No.~1, 163--171.
\par\noindent\hangindent=\parindent
11.~S.~A.~Kuleshov and D.~O.~Orlov, ``Exceptional sheaves over Del Pezzo
surfaces,'' {\it Russian Acad. Sci. Izv. Math}, {\bf44}(1995), No.~1, 479--513.
\par\noindent\hangindent=\parindent
12.~S.~A.~Kuleshov, ``Exceptional and rigid sheaves over surfaces with
anticanonical class free of basis components,'' Preprint No.~1 MK,
Independent University of Moscow, 1994 [in Russian].
\par\noindent\hangindent=\parindent
13.~S.~A.~Kuleshov, ``The new proof of the main theorem about
exceptional and rigid sheaves on $\PP^2$,'' Preprint MPI/95-11, 1995.
\par\noindent\hangindent=\parindent
14.~Yu.~I.~Manin, {\it Cubic Forms\/} North-Holland Math. Lib. (1974).
\par\noindent\hangindent=\parindent
15.~A.~A.~Markov, {\it On Binary Quadratic Forms of Positive Definition\/}
[in Russian], St.~Petersburg (1880).
\par\noindent\hangindent=\parindent
16.~D.~Yu.~Nogin, ``Helices of period four and equations of Markov type,''
{\it Math. USSR Izv., Ser.\ Math.}, {\bf37}, No.~1, 209--226 (1991).
\par\noindent\hangindent=\parindent
17.~D.~Yu.~Nogin, ``Helices on some Fano threefolds: Constructivity of
semiorthogonal bases of K$_{0}$,'' {\it Ann.\ ENS}, Ser.~4, {\bf27}, No.~2,
129--172 (1994).
\par\noindent\hangindent=\parindent
18.~A.~N.~Rudakov, ``The Markov numbers and  exceptional bundles on
${\Bbb P}^2$,'' {\it Math. USSR Izvestiya}, {\bf32}, (1989) No.~1, 99--112.
\par\noindent\hangindent=\parindent
19.~A.~N.~Rudakov, ``Exceptional bundles on a quadric,''
{\it Math. USSR Izvestiya}, {\bf33} (1989), No.~1, 115--138.

\end{document}